\documentclass[final,leqno]{siamltexmm2}
\usepackage{amsmath,amssymb}
\usepackage{graphicx}% Include figure files

\usepackage{algorithm}
\usepackage{algpseudocode}
\usepackage{tshortcuts}
\algnewcommand\algorithmicinput{\textbf{Input:}}
\algnewcommand\Input{\item[\algorithmicinput]}

\newcommand*\samethanks[1][\value{footnote}]{\footnotemark[#1]}

\begin{document}

%%%% Article title to be placed here
\title{Model selection for hybrid dynamical systems via sparse regression}

\author{%%%% Author details
N. M. Mangan$^{\dagger \S}$\thanks{Department of Engineering Sciences and Applied Mathematics, Northwestern University, Evanston, IL, 60208,  \email{niallmm@uw.edu}},
\and T. Askham\thanks{Department of Applied Mathematics, University of Washington, Seattle, WA. 98195},
\and S. L. Brunton\thanks{Department of Mechanical Engineering, University of Washington, Seattle, WA. 98195}, 
\and J. N. Kutz\samethanks[2] ,  
\and J. L. Proctor\thanks{Institute for Disease Modeling, Bellevue, WA 98005, USA }}

%%%%%%%%%%%%%%%%%%%%%%%%%%%

\maketitle
\newcommand{\slugmaster}{%
\slugger{siads}{xxxx}{xx}{x}{ }}%slugger should be set to juq, siads, sifin, or siims

%

%%%% Abstract text to be placed here %%%%%%%%%%%%
\begin{abstract}
Hybrid systems are traditionally difficult to identify and analyze using classical dynamical systems theory. 
Moreover, recently developed model identification methodologies largely focus on identifying a single set of governing equations solely from measurement data.
%
%For example, the {\sl sparse identification of nonlinear dynamics} (SINDy) methodology allows rapid model identification from a comprehensive library of models, but fails to identify hybrid systems and switching behaviors.  
%
In this article, we develop a new methodology, Hybrid-Sparse Identification of Nonlinear Dynamics (Hybrid-SINDy), which identifies separate nonlinear dynamical regimes, employs information theory to manage uncertainty, and characterizes switching behavior.
Specifically, we utilize the nonlinear geometry of data collected from a complex system to construct a set of coordinates based on measurement data and augmented variables.  
Clustering the data in these measurement-based coordinates enables the identification of nonlinear hybrid systems.  
This methodology broadly empowers nonlinear system identification without constraining the data locally in time and has direct connections to hybrid systems theory.  
We demonstrate the success of this method on numerical examples including a mass-spring hopping model and an infectious disease model. 
Characterizing complex systems that switch between dynamic behaviors is integral to overcoming modern challenges such as eradication of infectious diseases, the design of efficient legged robots, and the protection of cyber infrastructures.

%We believe automated, data-driven methods such as Hybrid-SINDy will provide intuition and predictive value to modern applications.  

\end{abstract}

%%%%%%%%%%%%%%% End of first page %%%%%%%%%%%%%%%%%%%%%

\section{Introduction}
\label{s:intro}
The high-fidelity characterization of complex systems is of paramount importance to manage modern infrastructure and improve lives around the world.
However, when a system exhibits nonlinear behavior and switches between dynamical regimes, as is the case for many large-scale engineered and human systems, model identification is a significant challenge.
%However, significant challenges face the identification of models representing large-scale, engineered and human systems that exhibit nonlinear behavior and switch between dynamical regimes.
%
These {\it hybrid systems} are found in a diverse set of applications including epidemiology~\cite{keeling2001seasonally}, legged locomotion~\cite{holmes2006dynamics}, cascading failures on the electrical grid~\cite{dobson2007complex}, and security for cyber infrastructure~\cite{li2014communication}.  
Typically, model selection procedures rely on physical principles and expert intuition to postulate a small set of candidate models; information theoretic approaches evaluate the goodness-of-fit to data amongst these models and penalizing over-fitting~\cite{Akaike1973,Akaike1974,nakamura2006comparative,lillacci2010parameter,penny2012comparing}.
Advances in data-driven methodologies, including the recently developed {\it sparse identification of nonlinear dynamics} (SINDy)~\cite{Brunton2016pnas}, broaden this procedure by selecting models from a combinatorially large library of possible nonlinear dynamical systems, decreasing the computational costs of model fitting and evaluations~\cite{Mangan2017}, and generalizing to a wide-variety of physical phenomena~\cite{mangan2016inferring,Rudy2017a}.  
However, standard model selection procedures and methods such as SINDy are not formulated to identify hybrid systems.
In this article, we describe a new method, called Hybrid-SINDy, which identifies hybrid dynamical systems, characterizes switching behaviors, and utilizes information theory to manage model selection uncertainty.

Predecessors of current data-driven model-selection techniques, called {\it system identification}, were developed by the controls community to discover linear dynamical systems directly from data~\cite{kalman:1965}.   
They made substantial advances in the model identification and control of aerospace structures~\cite{Juang1985jgcd,Phan1992jas}, and these techniques evolved into a standard set of engineering control tools~\cite{Katayama.2005}.  
One method to improve the prediction of linear input-output models was to augment present measurements with past measurements, i.e. delay embeddings \cite{Juang1985jgcd}.
Delay embeddings and their connections to Takens' embedding theorem have enabled equation-free techniques that distinguish chaotic attractors from measurement error in time series~\cite{Sugihara1990nature}, contribute to nonlinear forecasting~\cite{Sugihara1994ptrsla,brunton2017chaos}, and identify causal relationships among subsystems solely from time-series data~\cite{Sugihara2012}.  
Augmenting measurements with nonlinear transformations has also enabled identification of nonlinear dynamical systems from data.
As early as 1991, nonlinear feature augmentation was used to identify a nonlinear dynamical system for control~\cite{kowalski1991nonlinear}.  
Developed more recently, dynamic mode decomposition~\cite{Schmid2010jfm,Rowley2009jfm,kutz2016dynamic} has been connected to nonlinear dynamical systems via the Koopman operator~\cite{Mezic2005nd,Rowley2009jfm,Mezic2013arfm}.
More sophisticated data transformations, originating in the harmonic analysis community, are also being utilized for identifying nonlinear dynamical systems~\cite{COIFMAN20065,Giannakis2012pnas}.
Similarly, SINDy exploits these nonlinear transformations by building a library of nonlinear dynamic terms constructed using data.  This library is systematically refined to find a parsimonious dynamical model that represents the data with as few nonlinear terms as possible~\cite{Brunton2016pnas}.
%

%Methods like SINDy are currently not designed for hybrid systems, because they assume all time points come from a dynamical system with consistent equations. 
Methods like SINDy are currently not designed for hybrid systems, because they assume all measurement data in time is collected from a dynamical system with a consistent set of equations. 
In hybrid systems, the equations may change suddenly in time and one would like to identify the underlying equations without knowledge of the switching points.
%
%One approach is to constructing the models locally in time by {restricting} the input data to a short time-window.
One approach is to construct the models locally in time by {restricting} the input data to a short time-window.
Statistical models, such as the auto-regressive moving average (ARMA) and its nonlinear counterpart (NARMA), constrain the time-series to windows of data near the current time~\cite{box2015time}.
This technique has been extended to analyze non-autonomous dynamical systems, including hybrid systems, with Koopman operator theory~\cite{macevsic2017koopman}.
Alternatively, recent methods for recurring switching between dynamical systems use a Bayesian framework to infer how the state of the system, modeled as linear partitions, depends on multiple previous time-steps \cite{Linderman2016}.
This method enables reconstruction of state space in terms of linear generated states, and provides location-dependent behavioral states.

While restricting data locally in time may avoid erroneous model selection at the switching point, this method creates a new problem: there may not be enough data within a single window for data-driven model selection to robustly select and validate nonlinear models.
For nonlinear-model selection to work, we need to group together enough data from a consistent underlying model to perform selection and validation. 
Simplex-projection, which is used in cross convergent mapping, employs delay embeddings to find geometrically similar data for prediction \cite{Sugihara1990nature}.
Recently, Yair {\sl et al.} showed that data from dynamically similar systems could be grouped together in a label-free way by measuring geometric closeness in the data using a kernel method \cite{yair2017reconstruction}.
Here, we show that nonlinear model selection can succeed for hybrid dynamical systems when the data is examined within a pre-selected coordinate system that takes advantage of the intrinsic geometry of the data.

We present a generalization of SINDy, called Hybrid-SINDy, that allows for the identification of nonlinear hybrid dynamical systems.  
We utilize modern machine-learning methodologies to identify clusters within the measurement data augmented with features extracted from the measurements.  
Applying SINDy to these clusters generates a library of candidate nonlinear models.  
We demonstrate that this model library contains the different dynamical regimes of a hybrid system and use out-of-sample validation with information theory to identify switching behavior.
We perform an analysis of the effects of noise and cluster size on model recovery.
Hybrid-SINDy is applied to two realistic applications including legged locomotion and epidemiology.
These examples span two fundamental types of hybrid systems:  time and state dependent switching behaviors.

\section{Background}
\label{s:back}

\subsection{Hybrid Systems}
\label{ss:HybridBack}
Hybrid systems are ubiquitious in biological, physical, and engineering systems~\cite{keeling2001seasonally,holmes2006dynamics,dobson2007complex,li2014communication}.
Here, we consider hybrid models in which continuous-time vector fields describing the temporal evolution of the system state {\it change} at discrete times, also called events.  
Specifically, we choose a framework and definition for hybrid systems that is amenable to numerical simulations~\cite{back1993dynamical} and has been extensively adapted and utilized for the study of models~\cite{holmes2006dynamics}.
Note that these models are more complicated to define, numerically simulate, and analyze than classical dynamical systems with smooth vector fields~\cite{back1993dynamical,van2000introduction}.
Despite these challenges, solutions of these hybrid models have an intuitive interpretation: the solution is composed of piecewise continuous trajectories evolving according to vector fields that may change discontinuously at events.  
%Despite these challenges, solutions of these hybrid models have an intuitive interpretation: the solution is a concatenation of continuous trajectories evolving according to vector fields that are possibly different between events.  

Consider the state space of a hybrid system as a union 
\begin{equation}
V = \bigcup_{\alpha \in I}V_\alpha,
\end{equation}
where $V_\alpha$ is a connected open set in $\mathbb{R}^{n}$ called a chart and $I$ is a finite index.  
Describing the state of the system requires an index $\alpha$ and a point in $V_\alpha$, which we denote as $\mathbf{x}^\alpha$.  
We assume that the state within each patch evolves according to the classic description of a dynamical system $\dot{\mathbf{x}}^\alpha(t) = \mathbf{f}^\alpha(\mathbf{x}^\alpha(t))$, where $\mathbf{f}^\alpha(\mathbf{x}^\alpha)$ represents the governing equations of the system for chart $V_\alpha$.  
Transition maps $T^\alpha$ apply a change of states to boundary points within the chart; see \cite{holmes2006dynamics} for a more rigorous definition of $T^\alpha$.
%originally: changes of states
%
In this work, we consider hybrid systems where the transition between charts links the final state of the system on one chart $\mathbf{x}_{\alpha_i}$ to the initial condition on another $\mathbf{x}_{\alpha_j}$ where both $\mathbf{x}_{\alpha_i},\mathbf{x}_{\alpha_j} \in \mathbb{R}^n$.
Constructing the global evolution of the system {\it across} patches requires concatenating a set of smooth trajectories separated by a series of discrete events in time $\tau_1,\tau_2,\dots,\tau_o$. 
These discrete events can be triggered by either the state of the system $\tau_i(\mathbf{x})$ or external events in time $\tau_i(t)$.
In this article, we analyze hybrid systems representing both state-dependent and time-dependent events.
For a broader and more in-depth discussion on hybrid systems, we refer the reader to ~\cite{back1993dynamical,van2000introduction,holmes2006dynamics}.

\subsection{Sparse Identification of Nonlinear Dynamics (SINDy)}
\label{ss:SINDy}
%SINDy combines sparsity-promoting regression and nonlinear function library, to identify a dynamical system from time-series data. Here, we review the SINDy framework \cite{Brunton2016pnas}. The goal of SINDy is to find the functional form of dynamics, in terms of measured state variables $\mathbf{x} \in \mathcal{R}^n$. In the simplest case, we assume the equation takes the form of a few linearly additive functions, 
%\begin{equation}
%\frac{d}{dt} \mathbf{x}(t) = \sum\xi_l \mathbf{f}_l(\mathbf{x}(t))
%\label{eq:dynamical}
%\end{equation}
%where each $f_l(\mathbf{x})$ may be a nonlinear in terms of the state-variables. To identify appropriate functions from the data, we formulate a comprehensive library of candidate functions $\boldsymbol{\Theta}(\mathbf{x}) = [f_1(\mathbf{x}) \: f_2(\mathbf{x}) \: \dots \: f_p(\mathbf{x})] $, such that the true set are a subset of the columns in $\boldsymbol{\Theta}(\mathbf{x})$. The state variables are measured at $m$ time points, forming a time-series training-set $\mathbf{X} \in \mathcal{R}^{(m \times n)}$, where each row is a measurement of the state vector $\mathbf{x}^T(t_i)$, at time $t_i$. The function library is evaluated at the $m$ time points as $\boldsymbol{\Theta}(\mathbf{X}) \in \mathcal{R}^{(m \times p)}$. Derivative time-series data, $\dot{\mathbf{X}} \in \mathcal{R}^{(m \times n)}$, is measured or numerically calculated from $\mathbf{X}$. 
% In this subsection, we review the SINDy framework \cite{Brunton2016pnas}. 
%
SINDy combines sparsity-promoting regression and nonlinear function libraries to identify a nonlinear, dynamical system from time-series data \cite{Brunton2016pnas}.
We consider dynamical systems of the form
\begin{equation}
\frac{d}{dt} \mathbf{x}(t) = \sum_{l=1}^\zeta \xi_l \mathbf{f}_l(\mathbf{x}(t)),
\label{eq:dynamical}
\end{equation}
where $\mathbf{x}(t) \in \mathbb{R}^n$ is a vector denoting the state of the system at time $t$ and the sum of functions $\sum_{l=1}^\zeta \xi_l \mathbf{f}_l$ describes how the state evolves in time. 
Importantly, we assume $\zeta$ is small, indicating the dynamics can be represented by a parsimonious set of basis functions.
To identify these unknown functions from known measurements $\mathbf{x}(t)$, we first construct a comprehensive library of candidate functions $\boldsymbol{\Theta}(\mathbf{x}) = [\mathbf{f}_1(\mathbf{x}) \: \mathbf{f}_2(\mathbf{x}) \: \dots \: \mathbf{f}_p(\mathbf{x})] $.  
We assume that the functions in \eqref{eq:dynamical} are a subset of $\boldsymbol{\Theta}(\mathbf{x})$.
The measurements of the state variables are collected into a data matrix $\mathbf{X} \in \mathbb{R}^{(m \times n)}$, where each row is a measurement of the state vector $\mathbf{x}^T(t_i)$ for $i \in [1,m]$.
%
%The inter-measurement interval is equivalent for all measurements.
%
The function library is then evaluated for all measurements $\boldsymbol{\Theta}(\mathbf{X}) \in \mathbb{R}^{(m \times p)}$. 
The corresponding derivative time-series data, $\dot{\mathbf{X}} \in \mathbb{R}^{(m \times n)}$, is either directly measured or numerically calculated from $\mathbf{X}$. 

%We next pose a sparse-regression problem, which preserves the structure of Eqn \ref{eq:dynamical}:
%\begin{equation}
%\dot{\mathbf{X}} = \boldsymbol{\Theta}(\mathbf{X})\boldsymbol{\Xi},
%\end{equation} 
%where the unknown coefficients are in $\boldsymbol{\Xi}\in \mathcal{R}^{(p \times n)}$. Note that each column in $\boldsymbol{\Xi}$  determines the equation for a state variable, and we expect only a sparse set of the coefficients to be active in each column. Recovery of a sparse $\boldsymbol{\Xi}$ such that $||\dot{\mathbf{X}} - \boldsymbol{\mathbf{X}} \boldsymbol{\Xi}||_2$ is small, may be performed using sparse-regression methods, such as \textit{lasso}, elastic-net, or sequential least-squares with thresholding. 

To identify \eqref{eq:dynamical} from the data pair $(\boldsymbol{\Theta}(\mathbf{X}),\dot{\mathbf{X}})$, we solve
\begin{equation}
\dot{\mathbf{X}} = \boldsymbol{\Theta}(\mathbf{X})\boldsymbol{\Xi},
\end{equation} 
for the unknown coefficients $\boldsymbol{\Xi}\in \mathbb{R}^{(p \times n)}$ and enforce a penalty on the number of nonzero elements in $\boldsymbol{\Xi}$.     
Note that the $i$th column of $\boldsymbol{\Xi}$ determines the governing equation for the $i$th state variable.
We expect each coefficient vector in $\boldsymbol{\Xi}$ to be sparse, such that only a small number of elements are nonzero.
%
%Therefore, we are interested in solving the optimization problem 
%\begin{equation}
%\min_{\boldsymbol{\Xi}} \frac{1}{2} ||\dot{\mathbf{X}} - \boldsymbol{\Theta}(\mathbf{X})\boldsymbol{\Xi}||_2^2  \: \: \textrm{subject to} \: \: ||\boldsymbol{\Xi}||_0 \leq k,
%\label{eq:optprob}
%\end{equation}
%where the $l_0$ psuedo-norm of the coefficients, $||\boldsymbol{\Xi}||_0$, is the number of non-zero coefficients \cite{Bertsimas2016}. Finding the optimal subset of non-zero coefficients is computationally expensive for large $\boldsymbol{\Xi}$, and often a surrogate sparse-regression problem is solved instead. This has a general Lagrangian form  
We can find a sparse-coefficient vector using the Lagrangian minimization problem
\begin{equation}
\min_{\boldsymbol{\Xi}} \frac{1}{2} ||\dot{\mathbf{X}} - \boldsymbol{\Theta}(\mathbf{X})\boldsymbol{\Xi}||_2^2 + \hat{\lambda} R(\boldsymbol{\Xi}).
\label{eq:opt}
\end{equation}
Here, $R(\boldsymbol{\Xi})$ is a regularizing, sparse-penalty function in terms of the coefficients, and $\hat{\lambda}$ is a free parameter that controls the magnitude of the sparsity penalty. Two commonly used formulations include the \textit{LASSO} with an $l_1$ penalty $R(\boldsymbol{\Xi}) = ||\boldsymbol{\Xi}||_1$ and the elastic-net with an $l_1$ and $l_2$ penalty $R(\boldsymbol{\Xi}) = \gamma ||\boldsymbol{\Xi}||_1+ \frac{1}{2}(1-\gamma) ||\boldsymbol{\Xi}||_2^2 $ which includes a second free parameter $\gamma$ \cite{Hesterberg2008}. Less common, but perhaps more natural, is the choice $R(\bXi) = \|\bXi\|_0$, where the $\ell_0$ penalty is given by the number of nonzero entries in $\bXi$. In this article, we use sequential least-squares with hard thresholding to solve \eqref{eq:opt} with the $\ell_0$-type penalty, where any coefficients with values less than a threshold $\lambda$ are set to zero in each iteration \cite{Brunton2016pnas}. 

%Following this formulation of SINDy \cite{Brunton2016pnas}, several innovations have been made including extensions to partial differential equations \cite{Rudy2017a, Schaeffer2017a}, application to highly corrupted data \cite{Tran2016arxiv}, implicit formulation to handle rational functions \cite{Mangan2016ieee}, integral and weak formulations to improve robustness to noise \cite{Schaeffer2017b, Pantazis2017}, and use of group sparsity norms to recover structure in the presence of time-varying coefficients \cite{Schaeffer2017}. Additional connections with information criteria \cite{Mangan2016ieee}, and dynamic mode decomposition \cite{Kutz2016book} have also been explored. The connection with the Akaike information criteria (AIC) is essential for this work, as it allows automated evaluation of SINDy-generated models. 

Several innovations have followed the original formulation of SINDy~\cite{Brunton2016pnas}: the framework has been generalized to study partial differential equations~\cite{Rudy2017a, Schaeffer2017a} and systems with rational functional forms \cite{mangan2016inferring}; the impact of highly corrupted data has been analyzed~\cite{Tran2016arxiv}; the robustness of the algorithm to noise has been improved using integral and weak formulations~\cite{Schaeffer2017b, Pantazis2017}; and the theory has been generalized to non-autonomous dynamical system with time-varying coefficients using group sparsity norms~\cite{Schaeffer2017,Rudy2018arxiv}.  Additional connections with information criteria \cite{Mangan2017}, and extensions to incorporate known constraints, for example to enforce energy conservation in fluid flow models~\cite{Loiseau2017jfm}, have also been explored. The connection with the Akaike information criteria (AIC) is essential for this work, as it allows automated evaluation of SINDy-generated models.

\subsection{Model selection using AIC}
Information criteria provides a principled methodology to select between candidate models for systems without a well-known set of governing equations derived from first principles.
%. 
Historically, experts heuristically constructed a small number, $\mathcal{O}(10)$, of models based on their knowledge or intuition~\cite{burnham2002,Kuepfer2007,Hjorth2008,Schaber2011,woodward2004epidemiology,Blake22072014}. 
The number of candidate models is limited due to the computational complexity required in fitting each model, validating on out-of-sample data, and comparing across models.  
New methods, including SINDy, identify data-supported models from a much larger space of candidates without constructing and simulating every model~\cite{Schmidt2009science,Buchel2013,Cohen2015,Brunton2016pnas}. 
The fundamental goal of model selection is to find a parsimonious model, which minimizes error without adding unnecessary complexity through additional free parameters. 

In 1951, Kullback and Leibler (K-L) proposed a method for quantifying information loss or "divergence" between reality and model predictions \cite{Kullback1951}. 
Akaike subsequently calculated the relative information loss between models, by connecting K-L divergence theory with the likelihood theory from statistics. 
He discovered a deceptively simple estimator for computing the relative K-L divergence in terms of the maximized log-likelihood function for the data given a model, $L(\mathbf{x},{\hat{\mu}})$, and the number of free parameters, $k$~\cite{Akaike1973, Akaike1974}.
This relationship is now called Akaike's information criterion (AIC): 
\begin{equation}
AIC = 2 k  - 2 \ln ({L}(\mathbf{x},\hat{\mu})),
\label{Eq.AIC}
\end{equation}
where the observations are $\mathbf{x}$, and $\hat{ \mu}$ is the best-fit parameter values for the model given the data. The maximized log-likelihood calculation is closely related to the standard ordinary least squares when the error is assumed to be independently, identically, and normally distributed (IIND). In this special case, $AIC = \rho \ln(RSS/\rho) +2k$, where $RSS$ is the residual sum of the squares and $\rho$ is the number of observations.  The $RSS$ is expressed as $RSS = \sum_{i=1}^\rho (y_i- g(x_i;\mu))^2$ where $y_i$ are the observed outcomes, $x_i$ are the observed independent variables, and $g$ is the candidate model~\cite{burnham2002}. Note that the $RSS$ and the log-likelihood are closely connected.

In practice, the AIC requires a correction for finite sample sizes given by
\begin{equation}
AIC_c = AIC + 2(k+1)(k+2)/(\rho-k-2).
\label{Eq.AIC_rel}
\end{equation}
AIC and AIC$_c$ contain arbitrary constants that will depend on the sample size. These constants cancel out when the minimum AIC$_c$ across models is subtracted from the AIC$_c$ for each candidate model $j$, producing an interpretable model selection indicator called relative AIC$_c$, described by $\Delta{AIC}_c^j = AIC_c^j- AIC_c^{min}$. The model with the most support will have a score of zero; $\Delta{AIC}_c$ values allows us to rank the relative support of the other models.  Anderson and Burnham in their seminal work \cite{burnham2002} prescribe a general rule of thumb when comparing relative support among models:  models with $\Delta{AIC}_c <2$ have substantial support,  $4<\Delta{AIC}_c <7$ have some support, and $\Delta{AIC}_c >10$ have little support. 
These thresholds directly correspond to a standard p-value interpretation; we refer the reader to \cite{burnham2002} for more details.
In this article, we use $\Delta{AIC}_c=3$ as a slightly larger threshold for support in this study.
Following the development of AIC, many other information criteria have been developed including Bayesian information criterion (BIC)~\cite{schwarz1978estimating}, cross-validation (CV)~\cite{bishop2006pattern}, deviance information criterion (DIC)~\cite{linde2005dic}, and minimum description length (MDL)~\cite{rissanen1978modeling}. However, AIC remains a well-known and ubiquitous tool; in this article, we use relative AIC$_c$ with correction for low data-sampling~\cite{Mangan2017}.

\section{Hybrid-SINDy}
\begin{figure}[t]
	\centering
	\includegraphics[width=\textwidth]{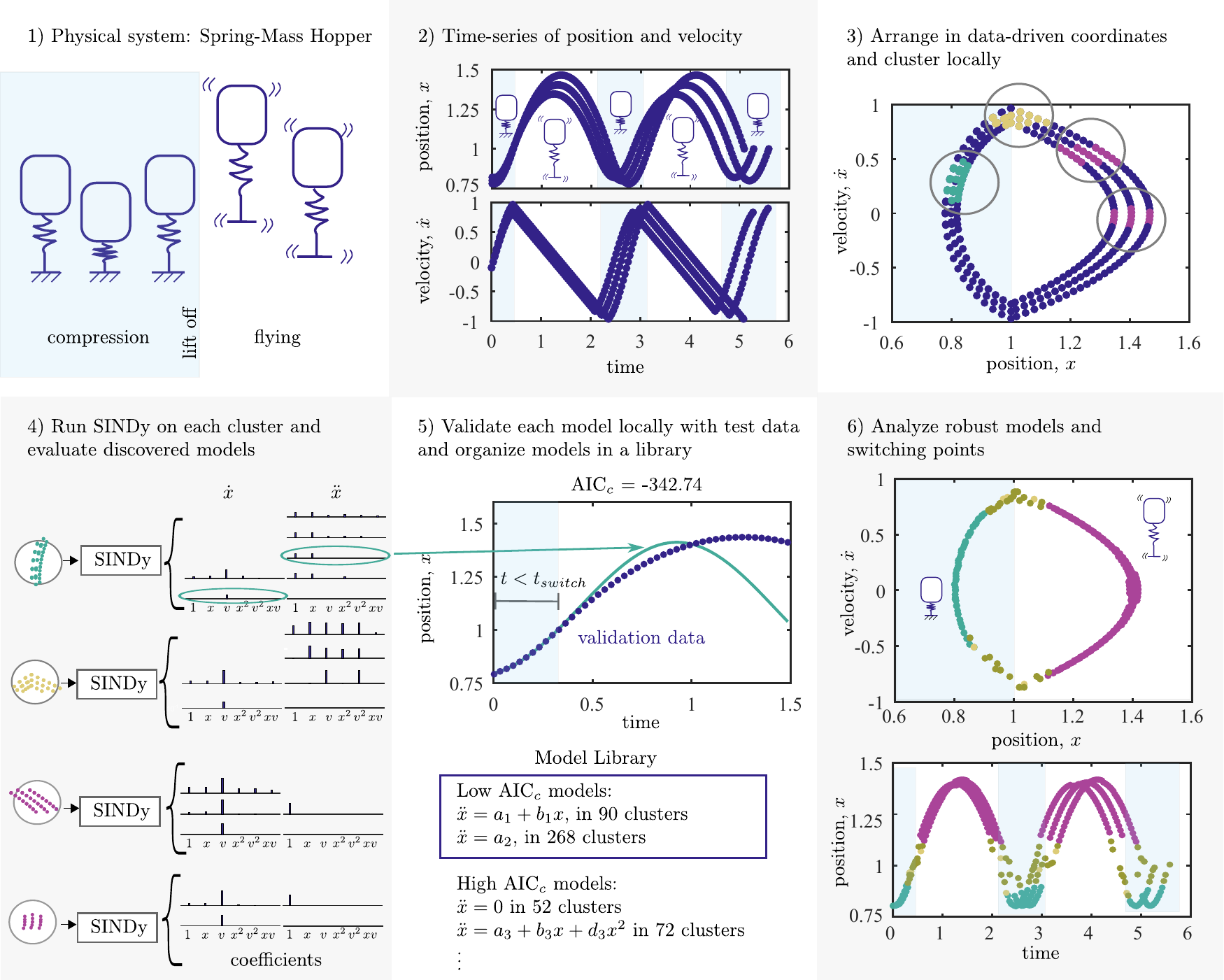}
	%%% where xxxxxx name represents "figurename.eps"
	\caption{Overview of the hybrid-SINDy method, demonstrated using the Spring-Mass Hopper system. Panel one illustrates the two dynamic regimes of spring compression (blue) and flying (white). Time-series for the position and velocity of the system sample both regimes (Panel 2). Clustering the data in data-driven coordinates allows separation of the regimes, except at transition points near $x = 1$ (Panel 3). Performing sparse model selection on each cluster produces a number of possible models per cluster (Panel 4). Panel 5 illustrates validating each model within the cluster to form a model library containing low AIC$_c$ models across all clusters. In panel 6, we plot the location of the 4 most frequent models across clusters. These models correctly identify the compression, flying, and transition points.}
	\label{fig:overview}
\end{figure}

Hybrid-SINDy is a procedure for augmenting the measurements, clustering the measurement and augmented variables, and selecting a model using SINDy for each cluster. We describe how to validate these models and identify switching between models. An overview of the hybrid-SINDy method is provided in Fig.~\ref{fig:overview} and Algorithm~\ref{alg}.

\begin{algorithm}[t]
	\caption{Hybrid-SINDy}
	\label{cSVDalgorithm}
	\begin{algorithmic}[1]
		\Input  The measurement data $\mathbf{X} \in \mathbb{R}^{b \times n}$, the set of measurement variables $D \in \mathbb{R}^{d \times 1}$ for clustering, the length of validation time-series $q$, the number of data points in the training set $m$, the number of data points in the validation set $v$, the sparsification values $\boldsymbol{\lambda} \in \mathbb{R}^r$, the number of library terms $p$, and the number of samples in each cluster $K$. 
		
		%		Training data and time derivative matrices $\mathbf{X}_T,\dot{\mathbf{X}}_T \in \mathbb{R}^{m \times n}$,  validation data and time derivative matrices $\mathbf{X}_V, \dot{\mathbf{X}}_V \in \mathbb{R}^{m \times n}$, the data-driven coordinates $D$ in terms of column indices for $\mathbf{X}_T = [\mathbf{X}_T~~ \dot{\mathbf{X}}_T]$, the input data in terms of the $D$ coordinates $\mathbf{Y}_T =\mathbf{X_T}(:,D)$ and $\mathbf{Y}_V =\mathbf{X_V}(:,D)$, and the number of points in the validation time-series $q$.
		%\Require $m \geq n$, integer $k \geq 1$ and $p \geq k$
		\Procedure{hybrid-SINDy}{$\mathbf{X}$,$D$,$s$,$m$,$v$,$c$,$K$}
		
		\State $\mathbf{X}_T\in \mathbb{R}^{m\times n},\mathbf{X}_V \in \mathbb{R}^{v\times n} \gets$ testTrainSeparation($\mathbf{X}$,$m$,$v$)
		\Comment{Construct training/validation}
		
		\State $\mathbf{\dot{X}}_T \in \mathbb{R}^{m\times n},\mathbf{\dot{X}}_V\in \mathbb{R}^{v\times n}	\gets$ derivative($\mathbf{X}_T,\mathbf{X}_V$)
		\Comment{Compute derivative matrix}
		
		\State $\mathbf{Y}_T \in \mathbb{R}^{m \times d}\gets$ variables($\mathbf{X}_T,\mathbf{\dot{X}}_T$,$D$)
		\Comment{Construct Augmented Measurements}
		\State $\mathbf{Y}_V \in \mathbb{R}^{m \times d} \gets$ variables($\mathbf{X}_V,\mathbf{\dot{X}}_V$,$D$)
		\Comment{Construct Augmented Measurements}
		
		%\For{$t_i \in \{t_1, t_2, \dots t_m\}$}
		\For{$i \in \{1, 2, \dots m\}$}
		\Comment{For each sample in the training set $t_i$, compute:}
		\State $\mathbf{C}_T^i \in \mathbb{R}^{K} \gets$ cluster($\mathbf{Y}_T$, $\mathbf{Y}_T(i,:)$,$K$)		
		\Comment{Cluster K samples from $\mathbf{Y}_T$ for each $\mathbf{Y}_T(i,:)$ }
		
		\State $\mathbf{C}_V^i \in \mathbb{R}^{K} \gets$ cluster($\mathbf{Y}_V$, centroid($\mathbf{Y}_T(\mathbf{C}_T^i,:)$),$K$)
		\Comment{Cluster K samples of $\mathbf{Y}_V$ for $\mathbf{C}_T^i$ }
		
		\State $\boldsymbol{\Theta}^i\in \mathbb{R}^{m \times p}$  $\gets$ library($\mathbf{X}_T(\mathbf{C}_T^i,:)$)
		\Comment{Generate library that contains $p$ features}
		
		%\For{$\lambda (j) \in \{ \lambda_1, \lambda_2, \cdots, \lambda_r \}$}  %Index over j instead??
		\For{$j \in \{ 1, 2, \cdots, r \}$}  %Index over j instead??
		\Comment{Search over sparsification parameter $\boldsymbol{\lambda}$.}
		
		\State Model($j$) $\gets$ SINDy($\mathbf{\dot{X}}_T(\mathbf{C}_T^i$,:), $\boldsymbol{\Theta}^i$, $\lambda(j)$ )
		\Comment{Identify sparse features \& model.}
		
		\For{$s \in \{1, 2, ..., K\}$}
		\Comment{Calculate error for each point in cluster}
		
		\State $\mathbf{Z}\in \mathbb{R}^{q \times n}$ $\gets$ simulate (Model($j$), $\mathbf{X}_V(\mathbf{C}_V^i(s),:),q)$
		\Comment{Simulate model}
		
		\State $ \mathbf{Z_V} \in \mathbb{R}^{q\times n} \gets \textnormal{find}(\mathbf{X}_V,\mathbf{X}_V(\mathbf{C}_V^i(s),:),q)$
		\Comment{Find validation time-series }

		\State $t_s \gets$ detect switching($\mathbf{Z},\mathbf{Z}_V$)
		\Comment{Find switching time}
		\For{$l \in \{1, 2, ..., n\}$}
		\Comment{Calculate error}
		\State		$E_{variable}(l) \gets$ $\frac{1}{t_s} \sum_{a=1}^{t_s} ( \mathbf{Z}(a,l)- \mathbf{Z}_V(a,l))^2$
		\Comment{Avg. over time}
		\EndFor
		\State $E_{avg}(s) \gets$ $\frac{1}{n} \sum_{l=1}^{n} E_{variable}(l)$ 
		\Comment{Avg. over measurements}
		\EndFor
		
		\State k$ \gets \text{numberOfFreeParameters}($Model$(j))$
		
		\State AIC$_{c}$($j$) $\gets$ ComputeAIC$_c$($E_{avg}$,k,K)
		%\Comment{Compute information criteria.}
		\EndFor
		\State $\Delta$ AIC$_c$ $\gets$ sort(AIC$_c$-minimum(AIC$_c$))
		\Comment{Rank models by relative AIC$_c$ scores.}
		\State  $\boldsymbol{\Pi} \gets$ $\mathbf{I}$(Model($\Delta$ AIC$_c$$<3$))
		\Comment{Store models with support in library}
		\EndFor
		
		\State ind $\gets$ sort(frequency($\boldsymbol{\Pi}$))
		\Comment Sort models by frequency across clusters.
		\State \textbf{return} $\boldsymbol{\Pi}(ind)$
		\Comment return the most frequent models.
		\EndProcedure
	\end{algorithmic}
	\label{alg}
\end{algorithm}

\subsection{Collect time-series data from system}
%Time-series points or samples from states, $\mathbf{x} \in \mathbb{R}^n$ of the dynamical system of interest are collected (Panel 2, Fig \ref{fig:overview}) from a system of interest (Panel 1, Fig \ref{fig:overview}). The training samples of the states are arranged into $\mathbf{X}_T \in \mathbb{R}^{(m\times n)}$ so that each row contains a measurement of the state variables in time, $\mathbf{x}(t_k)^T$. A column in $\mathbf{X}_T$ is a single state variable at each time sampled. We construct the testing or validation data $\mathbf{X}_V$ in a similar fashion. The derivatives of the training, $\mathbf{\dot{X}}_T$, and validation, $\mathbf{\dot{X}}_V$, data are either measured or calculated.

Discrete measurements of a dynamical system are collected and denoted by $\mathbf{x}(t_i) \in \mathbb{R}^n$; see Fig.~\ref{fig:overview}(2) for a time-series plot of the hopping robot illustrated in Fig.~\ref{fig:overview}(1). 
The measurement data is arranged into the matrix $\mathbf{X} =\left[ \mathbf{x}(t_1) ~~\mathbf{x}(t_2) ~\dots~ \mathbf{x}(t_{b}) \right]^T \in \mathbb{R}^{(b\times n)}$, where $^T$ is the matrix transpose.
The time-series may include trajectories from multiple initial conditions concatenated together.
The SINDy model is trained with a subset of the data $\mathbf{X}_T \in \mathbb{R}^{(m\times n)}$, where $m$ is the number of training samples.  
The corresponding data matrices for validation are denoted $\mathbf{X}_V \in \mathbb{R}^{v\times n}$, where $v$ is the number validation samples, and $b = m+v$.

\subsection{Clustering in measurement-based coordinates}

%The core insight of the method is to arrange the samples into a data-driven coordinate frame and find local clusters of samples in this space, as shown in Panel 3 of Fig. \ref{fig:overview}. By data-driven coordinates, we mean a geometric arrangement of the data such that samples from the same dynamical regimes are close to one another within the space. There are many existing algorithms for identifying such coordinates without {\sl a priori} knowledge of the number of dynamical regimes or which samples come from which regime \cite{Lee2007, Yair2016}. Here, the data-driven coordinates taken to be known as some subset of the state variables and their derivatives, indicated by the set of indices $D = {d} \in \mathbb{R}^{2n}$. These indices denote the columns in $\mathbf{Y}_T = [\mathbf{X}_T \; \dot{\mathbf{X}}_T],$ and $\mathbf{Y}_V = [\mathbf{X}_V \; \dot{\mathbf{X}}_V] \in \mathbb{R}^{m\times2n}$ which make up the data-driven coordinates.

%These measurements can be augmented with variables such as the derivative $\dot{\mathbf{x}}(t_i)$, nonlinear transformations $\mathbf{\tilde{f}}(\mathbf{x}(t_i))$~\cite{COIFMAN20065,Giannakis2012pnas}, or time-delay coordinates $[\mathbf{x}(t_i))^T~~ \mathbf{x}(t_{i-1}))^T~~ \dots]$~\cite{Juang1985jgcd,Sugihara1994ptrsla}.
%
%
%In principle, one would like to construct measurement-based coordinates utilizing the instrinsic, nonlinear geometry of the measurements.
%
Applications may require augmentation with variables such as the derivative, nonlinear transformations~\cite{COIFMAN20065,Giannakis2012pnas}, or time-delay coordinates~\cite{Juang1985jgcd,Sugihara1994ptrsla}.
In this article, we augment the state measurements  $\mathbf{x}(t_i)$ with the time derivative of the measurements.  
The time derivative matix is constructed similarly to the measurement matrix $\dot{\mathbf{X}}_T =\left[ \dot{\mathbf{x}}(t_1) ~~\dot{\mathbf{x}}(t_2) ~\dots~ \dot{\mathbf{x}}(t_m) \right]^T \in \mathbb{R}^{(m\times n)}$. 
The matrices $\dot{\mathbf{X}}_T$ and $\dot{\mathbf{X}}_V$ are either directly measured or calculated from $\mathbf{X}_T$ and $\mathbf{X}_V$, respectively.
If all state variables are accessible, such as in a numerical simulation, these data-driven coordinates directly correspond to the phase space of a dynamical system.  
Note that this coordinate system does not explicitly incorporate temporal information. 
Fig.~\ref{fig:overview}(3) illustrates the coordinates $(x,\dot{x})$ for the hopping robot.  
%
%Importantly, these coordinates are based on measurements of the system, but can be augmented by derivative information, delay embeddings~\cite{Juang1985jgcd,Sugihara1994ptrsla}, or more sophisticated transformations of the measurement data~.  
%
A subset of $[\mathbf{X}_T~~\dot{\mathbf{X}}_T]$ can also be utilized as measurement-based coordinates.  The set of indices $D$ are the measurements (columns), which are included in the analysis denoted by $\mathbf{Y}_T$ and $\mathbf{Y}_V$.

%We then construct clusters of samples in the training and validation sets. Starting with the first time point in the training set, $\mathbf{Y}_T(1,D)$ we use the nearest neighbor algorithm in MATLAB (knnsearch) to find a cluster of size $K$ in $\mathbf{Y}_T(i,D)$, excluding $i=1$. These training-set clusters $\mathbf{C}_T^i \in \mathbb{R}^{K \times n}$ are found for each time point $t_i \in t_1, t_2, \dots, t_m$. For each cluster in the training set, we calculate its centroid. We then find the cluster with $K$ samples, $\mathbf{C}_V^i \in \mathbb{R}^{K \times n}$, in the validation data set, $\mathbf{Y}_V(:,D)$ nearest the centroid of the training cluster. Panel 1 of Fig. \ref{fig:validation} illustrates a validation set in data-driven coordinates, with the centroids of 3 training clusters in black and the corresponding validation clusters in teal, gold, and purple. By finding these corresponding validation clusters, we can ensure our validation set is nearby. We construct a validation time-series from $\mathbf{Y}_V$ by selecting the $q$ time points starting from each measurement variable sample in $\mathbf{C}_V^i$, where $q<<m$. These short subsets of the validation time-series, stored in $\mathbf{X}_V^k \in \mathbb{R}^{q \times n}$ for the $k$th point in the cluster, are essential validation time-series that have initial conditions within cluster $i$. 

We then identify clusters of samples in the training and validation sets.  For each sample (row) in $\mathbf{Y}_T$, we use the nearest neighbor algorithm {\it knnsearch} in MATLAB to find a cluster of $K$ similar measurements in $\mathbf{Y}_T$.  
These training-set clusters, which are row indices of $\mathbf{Y}_T$ denoted $\mathbf{C}_T^i \in \mathbb{R}^{K}$, are found for each time point $t_i \in \left[ t_1, t_2, \dots, t_m \right ]$.  The centroid of each cluster is computed within the training set $\mathbf{Y}_T(\mathbf{C}_T^i)$. 
We then identify $K$ measurements from $\mathbf{Y}_V$ near the {\bf training} centroid clusters.
Note that these clusters in the validation data, $\mathbf{C}_V^i \in \mathbb{R}^{K}$, are essential to testing the out-of-sample prediction of Hybrid-SINDy.
Fig. \ref{fig:validation}(1\&2) illustrates the validation set in measurement-based coordinates, with the centroids of 3 training clusters as black dots and the corresponding validation clusters in teal, gold, and purple dots. 

By finding the corresponding validation clusters, we ensure the out-of-sample data for validating the model has the same local, nonlinear characteristics of the training data.  To assess the performance of the models, we also need to identify a validation time-series from $\mathbf{Y}_V$.  Starting with each data point in a validation cluster $\mathbf{C}_V^i$, we collect $q$ measurements from  $\mathbf{Y}_V$ that are temporally sequential, where $q\ll m$.  These subsets of validation time-series, $\mathbf{Z}_V^i \in \mathbb{R}^{q \times n}$ are defined for each data point and each cluster.  The validation time-series helps characterize the out-of-sample performance of the model fit.

\subsection{SINDy for clustered data}
%We perform SINDy on each training cluster $\mathbf{C}_T^i$, using an alternating least squares and hard thresholding as in \cite{Brunton2016pnas}. For each cluster, we search over the sparsification parameter, $\lambda(j) \in \{ \lambda_0, \lambda_1, \dots \lambda_q\}$, generating many possible models for $\mathbf{\dot{x}}(t)$ per cluster (Panel 4, Fig. \ref{fig:overview}), as illustrated in Fig. \ref{fig:overview} Panel 4. There are often redundant models for multiple values of $\lambda$, so the number of models per cluster is generally less than $q$. 

We perform SINDy for each training cluster $\mathbf{C}_T^i$, using an alternating least squares and hard thresholding described in \cite{Brunton2016pnas} and \S\ref{s:back}\ref{ss:SINDy}. For each cluster, we search over the sparsification parameter, $\lambda(j) \in \{ \lambda_1, \lambda_1, \dots \lambda_r\}$, generating a set of candidate models for each cluster; see Fig.~\ref{fig:overview}(4) for an illustration.  In practice, the number of models per cluster is generally less than $r$ since multiple values of $\lambda$ can produce the same model.
In this article, the library, $\boldsymbol{\Theta}(\mathbf{X})$, includes polynomial functions of increasing order (i.e., $x,x^2,x^3,\dots$), similar to the examples in~\cite{Brunton2016pnas}.  
However, the SINDy library can be constructed with other functional forms that reflect intuition about the underlying process and measurement data.  

\subsection{Model validation and library construction}
%We perform a validation step on each model in the cluster and build a library models identified across all clusters. Using the validation cluster as a set of initial conditions, $\mathbf{C}_V^i$, we simulate each potential model in cluster $i$ for $q$ time steps producing simulated time-series $\mathbf{X}'$. We would like to compare these simulations against the local validation time-series $\mathbf{X}_V^k$ and calculate a AIC$_c$ score. A plot of  $\mathbf{X}_V^k$ and $\mathbf{X}'$ for a single cluster is shown in Panel 1 of Fig. \ref{fig:validation}. 

\begin{figure}[t]
	\centering
	\includegraphics[width=\textwidth]{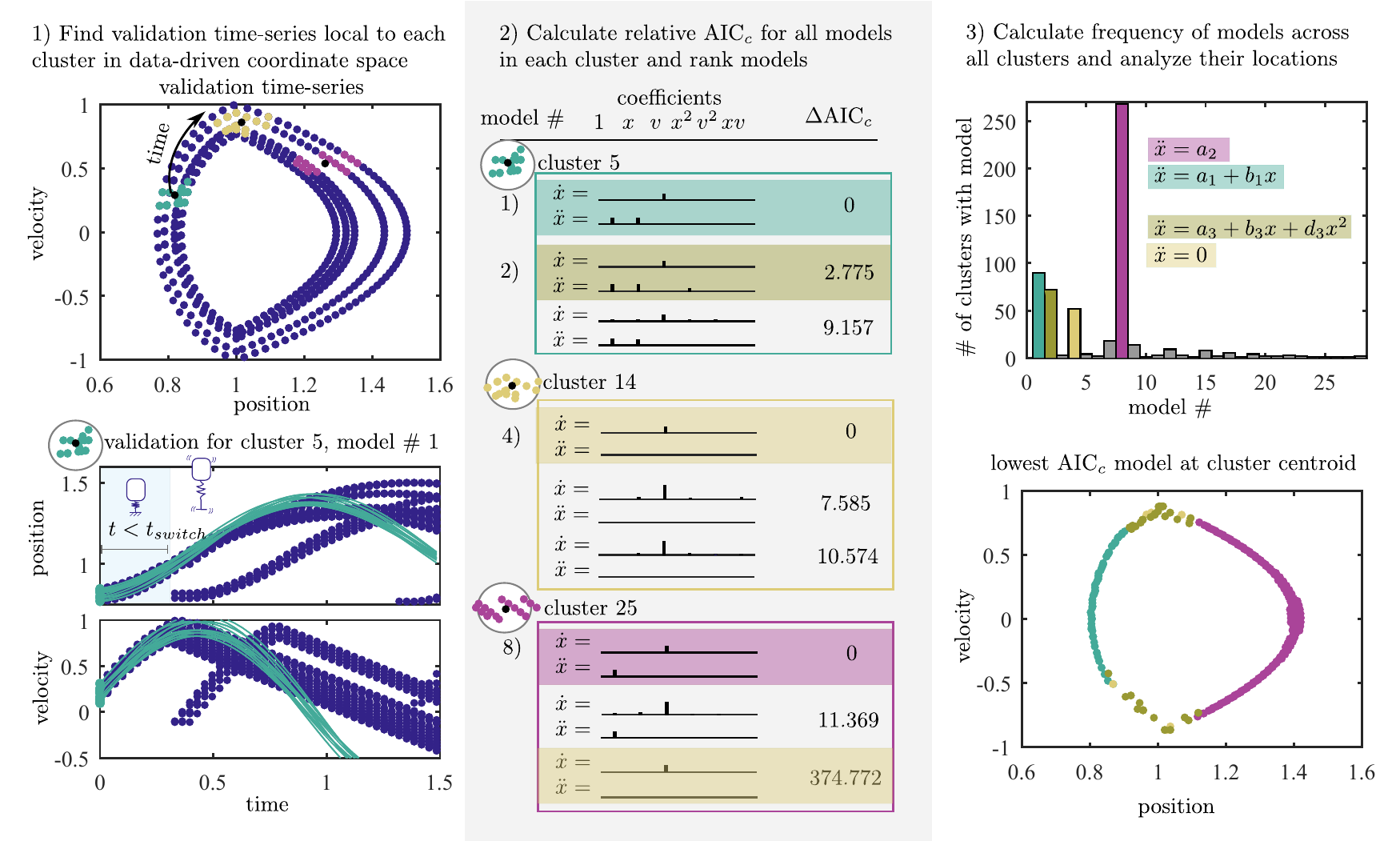}
	%%% where xxxxxx name represents "figurename.eps"
	\caption{Steps for local validation and selection of models. For each cluster from the training set, we identify validation time-series points that are local to the training cluster centroid (black dots, panel 1). We simulate time-series for each model in the cluster library, starting from each point in the validation cluster (teal, gold and purple dots) and calculate the error from the validation time-series. Using this error we calculate a relative AIC$_c$ value and rank each model in the cluster (panel 2). We collect the models with significant support into a library, keeping track of their frequency across clusters. The highest frequency models across clusters are shown in panel 3.  Note that the colors associated with each model in panel 3 are consistent across panels. }
	\label{fig:validation}
\end{figure}

%In this subsection, we describe the validation procedure for models within each cluster and the construction of a library of trusted models across all clusters.  Generally, the 
Validation involves producing simulations from candidate models and comparing to the validation data.  
Using the validation cluster as a set of $K$ initial conditions $\mathbf{C}_V^i$, we simulate each candidate model $j$ in cluster $i$ for $q$ time steps producing time-series $\mathbf{Z} \in \mathbb{R}^{q \times n}$.
We compare these simulations against the validation time-series $\mathbf{Z}_V$ and calculate an out-of-sample AIC$_c$ score. An example illustration comparing $\mathbf{Z}_V$ and $\mathbf{Z}$ for a single cluster is shown in Fig.~\ref{fig:validation}(1). 

%In order to calculate the error between the simulation and validation, we must account for the fact that the underlying dynamics may switch somewhere in the middle of the validation set.   We cannot know when this will happen {\sl a priori}.  To calculate a lower bound on the duration of time we can compare the simulation and validation data, $t_s$, we use an algorithm which detects changes in the mean of the absolute error between the simulated time-series and validation time-series. We use the findchangepoints function in MATLAB \cite{Killick2012}. Notably, this algorithm does not robustly find the time at which our time-series switch dynamical regimes. It will often prematurely identify a transition, especially in oscillatory systems. For this reason, we use it as a {\sl lower bound} where we can reasonably compare the simulated and validation data. 

In order to calculate the error between the simulation and validation, we must first account for the possibility of the dynamics switching before the end of the $q$ validation time-steps.  
We use the function {\it findchangepoints} in MATLAB~\cite{killick2012} to detect a change in the mean of the absolute error between the simulated and validation time-series.  The time index closest to this change is denoted $t_s$.
Notably, this algorithm does not robustly find the time at which our time-series switch dynamical regimes.  The algorithm tends to identify the transition prematurely, especially in oscillatory systems.  We utilize $t_s$ as a {\sl lower bound}, before which we can reasonably compare the simulated and validation data. 

%Once the subset of the simulated, $\mathbf{X}'(1:t_s,:)$, and validation, $\mathbf{X}_V^k(1:t_s,:)$, time-series have been identified for a particular model we calculate the average absolute error $\epsilon(k) = \frac{1}{t_s}\Sigma_0^{t_s}| \mathbf{X}'(1:t_s,:)- \mathbf{X}_V^k(1:t_s,:)|$. We calculate this average error for all $K$ initial conditions in the validation cluster, where each time series $\mathbf{X}_V^k$ may have its own switching point. Given these error values and the number of terms in the $j$th model we calculate an AIC$_c$ value \cite{Akaike1973,Akaike1974}. This is following the same method in \cite{Mangan2017}. We repeat this AIC$_c$ calculation for each model. Once we have AIC$_c$ scores for each model within the cluster, we calculate the relative AIC$_c$ scores and identify models within the cluster with significant support or relative AIC$_c < 3$ as shown in Panel 2 of Fig. \ref{fig:validation}. These models are used to build a model library and higher AIC$_c$ models are discarded (Panel 5, Fig. \ref{fig:overview}). The model library records the structure of selected models and how many times they appear as highly supported across clusters. 

To assess a model's predictive performance within a cluster, we compare the simulated data $\mathbf{Z}$ and validation data $\mathbf{Z}_V$ restricted to time points before $t_s$.  
Specifically, we calculate the residual sum of square error for a candidate model by comparing the $K$ time-series from the validation data to the model outputs, described by $E_{avg}(s) =\frac{1}{n}\sum_{l=1}^n \left (\frac{1}{t_s}\sum_{a=1}^{t_s}( z_{a,l}- z_{a,l}^V)^2 \right )$ for $s \in [1,K]$, where $z_{a,l}$ corresponds to the $a$ row and $l$ column of $\mathbf{Z}$ and similarly with $ z_{a,l}^V$ to $\mathbf{Z}_V$.
Thus, the vector $E_{avg}$ contains the average error over time points and state variables for $K$ initial conditions of model $r$.  
%
%Note that each time series in $\mathbf{Z}_V$ may have its own switching point.
%
For each candidate model $r$, we calculate the AIC$_\text{c}$ from \eqref{Eq.AIC_rel} using AIC$(r) = K \ln (E_{avg}/K) + 2 k$, the number of initial conditions in the validation set $K$, the average error for each initial condition, $E_{avg}$, and the number of free parameters (or terms) in the selected model $k$~\cite{Akaike1973,Akaike1974}.  
An equivalent procedure is found in \cite{Mangan2017}.  Once we have AIC$_c$ scores for each model within the cluster, we calculate the relative AIC$_c$ scores and identify models within the cluster with significant support where the relative AIC$_c < 3$; see Fig. \ref{fig:validation}(2) for an illustration. These models are used to build the model library. Models with larger relative AIC$_c$ are discarded, illustrated in Fig.~\ref{fig:overview}(5). The model library records the structure of highly-supported models and how many times they appear across clusters. 

\subsection{Identification of high-frequency models and switching events}
After building a library of strongly supported models, we analyze the frequency of model structures appearing across clusters, illustrated in Fig.~\ref{fig:validation}(3). 
The most frequent models and the location of their centroids provide insight into connected regions of measurement space with the same model; see Fig. \ref{fig:overview}(6) for an example.
By examining the location and absolute AIC$_c$ scores of the models, we can identify regions of similar dynamic behavior and characterize events corresponding to dynamic transitions.

\section{Results: Model selection }
\label{s:results}
\subsection{Mass-Spring Hopping Model}
\label{ss:massspring}
In this subsection, we demonstrate the effectiveness of Hybrid-SINDy by identifying the dynamical regimes of a canonical hybrid dynamical system: the spring-mass hopper.  The switching between the flight and compression stages of the hopper depends on the state of the system~\cite{holmes2006dynamics}.  Fig.~\ref{fig:overview}(2) illustrates the flight and spring compression regimes and dynamic transitions.  Note these distinct dynamical regimes are called charts, and liftoff and touchdown points are state-dependent events separating the dynamical regimes; see \S\ref{s:back}\ref{ss:HybridBack} for connections to hybrid dynamical system theory.  
The legged locomotion community has been focused on understanding hybrid models due to their unique dynamic stability properties~\cite{ruina1998nonholonomic}, the insight into animal and insect locomotion~\cite{full1999templates,holmes2006dynamics}, and guidance on the construction and control of legged robots~\cite{saranli2001rhex,raibert2008bigdog}.  

%Our first test case for the hybrid-SINDy algorithm is a spring-mass hopper system. We have selected the spring-mass hopper system as an example, because it is one of the canonical hybrid-dynamical systems and has been rigorously As shown in \ref{fig:overview}, this dynamical system is composed of two dynamical regimes: compression and flying. The simplest model of this system models the forces of gravity and a linear spring term during compression.
A minimal model of the spring-mass hopper is given by the following: 
\begin{equation}
m \ddot{x}  = 
\begin{cases}
-k(x-x_0) -mg, \:  & x \leq x_0  \\
-mg, \:  & x>x_0 
\end{cases}
\label{eq:hopper}
\end{equation}
where $m$ is the mass, $k$ is the spring constant, and $g$ is gravity. The unstretched spring length $x_0$ defines the flight and compression stages, i.e., $x>x_0$ and $x\leq x_0$, respectively.  For convenience, we non-dimensionalize \eqref{eq:hopper} by scaling the height of the hopper by $y = \frac{x}{x_0}$, scaling time by $\tau = t \sqrt{(\frac{k x_0}{m})}$, and forming the non-dimensional parameter $\kappa = \frac{k x_0}{m g}$.  Thus, $\kappa$ represents the balance between the spring and gravity forces.  Eq. \eqref{eq:hopper} becomes 
\begin{equation}
\ddot{y} =
\begin{cases}
1- \kappa(y-1), \; & y\leq 1  \\ 
-1, \; & y>1. 
\end{cases}
\label{eq:hopper2}
\end{equation} 
For our simulations we chose $\kappa = 10$. The switching point between compression and flying occurs at $y =1$ in this non-dimensional formulation. 

%The system is under compression when the height of the hopper, with mass $m$, is below the critical or equilibrium height, $x_0$, and flying when above. The hookian-spring constant for the system is $k$, and force of gravity is $g$. We can non-dimensionalize this system by scaling the height of the hopper by $y = \frac{x}{x_0}$, scaling time by $\tau = t \sqrt(\frac{k x_0}{m})$ and forming the non-dimensional parameter $\kappa = \frac{k x_0}{m g}$, representing the balance between the force of the spring and the force of gravity. With these substitutions, Eqn. \ref{eq:hopper} become $\ddot{y} = 1- \kappa*(y-1)$ and $\ddot{y} = -1$. For our simulations we chose $\kappa = 10$. The non-demensionalization means that the switching point between compression and flying occurs at $y =1$. 

\subsubsection{Generating input time-series from the model}

We generate time-series samples from \eqref{eq:hopper2} by selecting three initial conditions $(y_0,\dot{y}_0) \in \{(0.8, -0.1), (0.78, -0.1), (0.82,-0.1)\}$. 
We simulate the system for a duration of $t = [0~~5]$ with sampling intervals of $\Delta \tau =0.033$, producing 152 samples per initial condition. 
The resulting time-series of the position and velocity, $\mathbf{y}(t_i)$ and $\dot{\mathbf{y}}(t_i)$, are used to construct the training set matrices $\mathbf{X}_T = \mathbf{Y}_T = \left [ \begin{array}{cc}\mathbf{y}(t_i) & \dot{\mathbf{y}}(t_i) \end{array}\right]$  where each row corresponds to sample.
The position and velocity time-series are plotted in Fig.~\ref{fig:overview}(2).  Fig. \ref{fig:overview}(3) illustrates the position-velocity trajectories in phase-space.  We also add gaussian noise with mean zero and standard deviation $10^{-6}$ to the position and velocity time-series in $\mathbf{Y}_T$.  In this example, the derivatives $\dot{\mathbf{X}}_T = \dot{\mathbf{Y}}_T = [\dot{\mathbf{y}}(t_i)^T \; \ddot{\mathbf{y}}(t_i)^T]$ are computed exactly, without noise. The validation set $\mathbf{Y}_V$ is generated using the same intervals and duration, but for initial conditions: $(y_0,\dot{y}_0) \in \{(0.84, -0.11), (0.77, -0.12), (0.83,-0.13), (0.79, -0.13), (0.79, -0.10), (0.82,-0.11)\}$.

%To generate a time-series suitable for equation discovery, we chose three initial conditions, close together in position and velocity. The resulting time-series are plotted in Fig. \ref{fig:overview} step 2, and their phase-space representation in Fig. \ref{fig:overview} step 3. We add gaussian noise with mean zero and standard deviation $10^{-6}$ to the position and velocity time-series for the SINDy feature library. For this proof of concept, the derivatives used in the regression are exact, and not calculated from the noisy data. 

\subsubsection{Hybrid-SINDy discovers flight and hopping regimes}

In this case, the position and velocity measurements in phase space provide a natural, data-driven coordinate system to cluster samples.  Here, we identify $m = 492$ clusters, one for each timepoint.
% In this case, the indices of the training set, $\mathbf{Y}_T$, defining the data-driven coordinate system are $D = \{1,2\}$. 
Fig. \ref{fig:validation}(1) illustrates three of these clusters. We use a model library containing polynomials up to 2nd order in terms of $\mathbf{X}_T$.  Applying SINDy to each cluster, we produce a set of models for each cluster and rank them within the cluster using relative AIC$_c$; this procedure is illustrated in Fig. \ref{fig:validation}(2).  We retain only the models with strong support, relative AIC$_c < 3$.  Fig. \ref{fig:validation}(3) shows that the correct models are the most frequently identified by Hybrid-SINDy.  In addition, when we plot the location of the discovered models in data-driven coordinates (phase space for this example), we clearly identify the compression model when $y<1$ (teal) and the flying model when $y >1$ (purple). There is a transition region at $y=1$, where the incorrect models, plotted in gold and yellow, are the lowest AIC$_c$ models in the cluster. 

%From the phase-space or 'data-driven coordinate' representation, one can see that the choice of initial conditions allows for reasonable clusters. Applying SINDy to each cluster and validating the resulting models on validation time-series, shown in Fig. \ref{fig:validation} step 1, we generate a set of models for each cluster and rank them within the cluster using AIC$_c$, as shown in Fig. \ref{fig:validation} step 2. Taking the lowest AIC$_c$ models within each cluster, Fig. \ref{fig:validation} step 3 shows that the two correct models are the most frequently occurring across all the models. In addition, when we plot the location of these discovered model in data-driven coordinates we clearly identify the compression model is found when $y<1$ (teal) and the flying model is found when $y >1$ (purple). There is a transition region at $y=1$, where the incorrect models, plotted in gold and yellow, are the lowest AIC$_c$ models in the cluster. 

\begin{figure}[t]
	\centering
	\includegraphics[width=\textwidth]{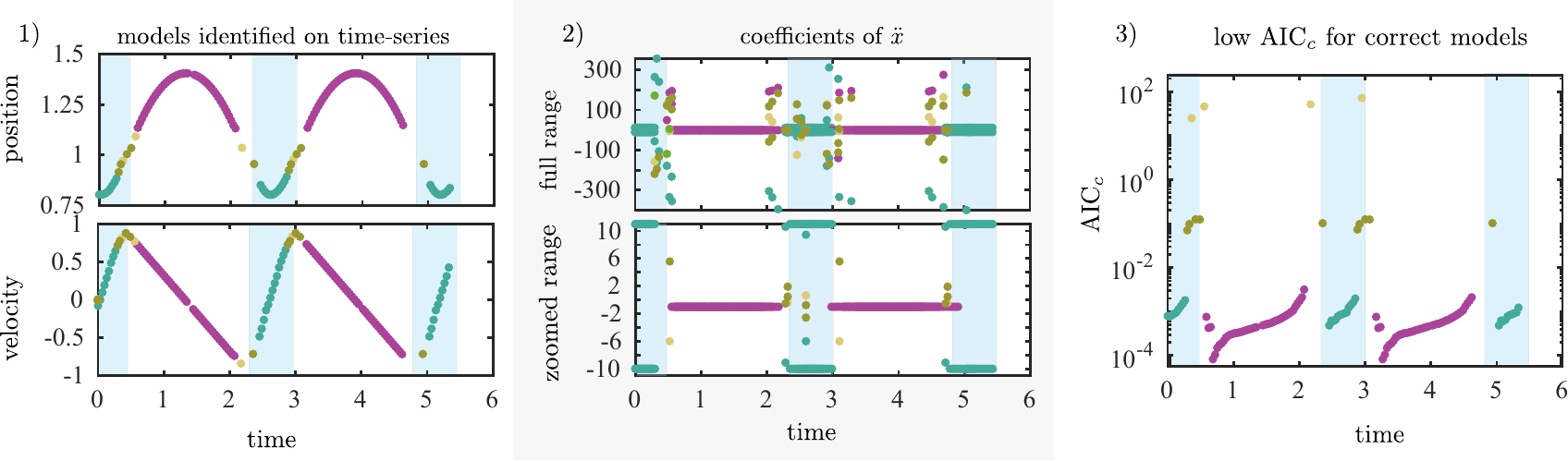}
	%%% where xxxxxx name represents "figurename.eps"
	\vspace*{-.1in}
	\caption{Hopping model discovery shown in time. Single time-series of data and associated model coefficients and AIC$_c$ are plotted as function of time with the correct models indicated by color. Teal dots indicate the recovery of compressed spring model, purple dots indicate recovery of the flying model and yellow and gold dots indicate recovery of incorrect models. Both coefficients and absolute AIC$_c$ are plotted for the model with only the lowest AIC$_c$ value at each cluster.}
	\label{fig:HoppingRes}
\end{figure}

To investigate the success of model discovery over time, Fig.~\ref{fig:HoppingRes} illustrates the discovered models (same color scheme as in Fig.~\ref{fig:validation}), the estimated model coefficients, and the associated absolute AIC$_c$ values. The four switching points between compression (teal) and flying (purple) area clearly visible, with incorrect models (gold and yellow) marking each transition.The model coefficients are consistent within either the compression or flying region, but become large within the transition regions, seen in Fig.~\ref{fig:HoppingRes}(2). 

The AIC$_c$ plot shows only the lowest absolute AIC$_c$ found in each cluster for the top 4 most frequent models across clusters. There is a substantial difference between the AIC$_c$ values for the correct (AIC$_c$ $\leq 3\times 10^{-3}$) and incorrect models (AIC$_c$ $\geq 2\times 10^{-2}$). As the system approaches a transition event, the AIC$_c$ for Hybrid-SINDy increases significantly. The increase is likely due to two factors: 1) as we approach the event there are fewer time points contributing to the AIC$_c$ calculation, and 2) we find an inaccurate proposed switching point, $t_s$, for validation data. As a switching event is approached, locating $t_s$ becomes challenging, and points from after a transition are occasionally included in the local error approximation $\epsilon(k)$. Note that the increase in AIC$_c$ between clusters provides a more robust indication of the switch than MATLAB's built in function {\it findchangepoints} applied to the time-series without clustering.  The {\it findchangepoints}, which uses statistical methods to detect change points, often fails for dynamic behavior such as oscillations.

\subsection{SIR disease model with switching transmission rates}

\begin{figure}[t]
	\centering
	\includegraphics[width=\textwidth]{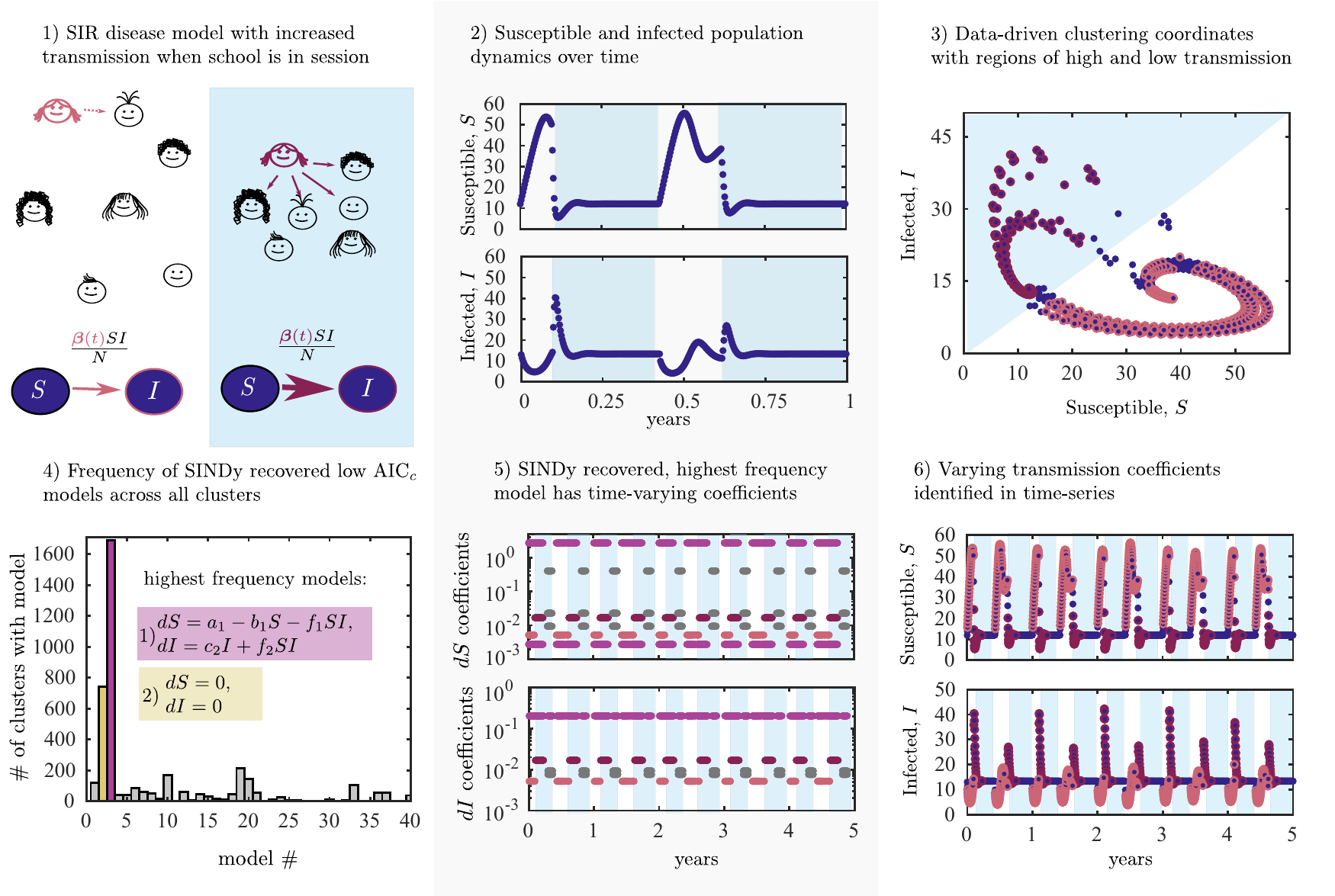}
	%%% where xxxxxx name represents "figurename.eps"
	\vspace*{-.1in}
	\caption{Sparse selection of Susceptible-Infected-Recovered (SIR) disease model with varying transmission rates. Panel 1: School children have lower transmission rates during school breaks (white background), and higher transmission due to increased contact between children while school is in session (blue background). The infected, $I$, and susceptible, $S$, population dynamics over one school year, show declines in the infected population while school is out of session, followed by spikes or outbreaks when school is in session as shown in panel 2. Clustering in data-driven coordinates $S$ vs $I$, shown in panel 3, and performing SINDy on the clusters, identifies a region with high transmission rate (maroon) and low transmission rate (pink). A frequency analysis across all clusters of the low AIC$_c$ models in each cluster, shown in panel 4, identifies 2 models of interest. The highest frequency model is the correct model, and SINDy has recovered the true coefficients for this model in both high and low transmission regimes. Panel 5 shows the coefficients of the highest frequency model recovered in time. Panel 6 overlays the recovered transmission rates on the time-series data used for selection.}
	\label{fig:SIR}
\end{figure}

In this subsection, we investigate a time-dependent hybrid dynamical system.  Specifically, we focus on the Susceptible, Infected and Recovered (SIR) disease model with varying transmission rates.  This dynamical system has been widely studied in the epidemiological community due to the nonlinear dynamics~\cite{keeling2001seasonally} and the related observations from data~\cite{grenfell1994measles}.  For example, the canonical SIR model can be modified to increase transmission rates among children when school is in session due to the increased contact rate~\cite{broadfootmeasles}. Fig. \ref{fig:SIR}(1) illustrates the switching behavior.  The following is a description of this model:
%For our second example, we model the Susceptible, Infected and Recovered (SIR) disease model with varying transmission rates. Temporal variation in transmission rates of diseases is a known issue in epidemiology [Citation]. A classic example is the increase in transmission rate through populations of children when school is in session [Citation]. Fig. \ref{fig:SIR} step 1 shows a schematic of this model. When children are out of school they have less contact with one another and transmission is slow. During the months children are in-school, the transmission rate is higher because contact rate between individuals increases as shown in Fig. \ref{fig:SIR} section 1. We simulate this system using the simple continuous ODE model for $S$, the susceptible population, $I$ the infected population, and $R$, the recovered population;
%

\begin{subequations}
	\begin{align}
	\dot{S} & =  \nu N - \frac{\beta(t)}{N} I S - d S  \label{eq:S} \\
	\dot{I} & = \frac{\beta(t)}{N} I S - (\gamma + d) I \label{eq:I} \\ 
	\dot{R} & = \gamma I - d R,
	\end{align}
	\label{eq:R}
\end{subequations}

where $\nu = \frac{1}{365}$ is the rate which students enter the population, $d = \nu$ is the rate at which students leave the population, $N = 1000$ is the total population of students, and $\gamma = \frac{1}{5}$ is the recovery rate when 5 days is the average infectious period.  The time-varying rate of transmission, $\beta(t)$, takes on two discrete values when school is in or out of session:
\begin{equation}
\beta(t) = 
\begin{cases}
\hat{\beta}(1+b), & \text{ \space t $\in$ school in session } \\
\hat{\beta}\frac{1}{(1+b)},  & \text{ \space t $\in$ school out of session. } 
\end{cases}
\end{equation}
The variable $\hat{\beta} = 9.336$ sets a base transmission rate for students and $b = 0.8$ controls the change in transmission rate. The school year is composed of in-class sessions and breaks.  The timing of these periods are outlined in Table \ref{tab:calendar}.  We chose these slightly irregular time periods, creating a time series with annual periodicity, but no sub-annual periodicity. A lack of sub-annual periodicity could make dynamic switching hard to detect using a frequency analysis alone.

\begin{table}[!h]
	\caption{School calendar for a year.} \label{tab:calendar}
	\center
	%	\vspace{-3mm}	
	\begin{tabular}{|l|c|c|c|}
		\hline
		session & days & time period & transmission rate \\
		\hline 
		winter break  & 0 to 35  & 1.2 months  & $\beta = 5.2$ \\
		\hline 
		spring term  & 35 to 155  & 4 months & $\beta = 16.8$ \\
		\hline
		summer break  & 155 to 225  & 2.3 months  & $\beta = 5.2$ \\
		\hline 
		fall term  & 225 to 365  & 4.6 months & $\beta = 16.8$ \\
		\hline
	\end{tabular}
\end{table}

\subsubsection{Generating input time-series from the SIR model}

To produce training time-series, we simulate the model for 5 years, recording at a daily interval.  This produces 1825 time-points.  We collect data along a single trajectory starting from the initial condition at $S_0 =12$ $I_0 =13$, $R_0= 975$. For this model, the dynamic trajectory rapidly settles into a periodic behavior, where the size of spring and fall outbreaks are the same each year.  
We add a random perturbation to the start of each session by changing the number of children within the the $S$, $I$, and $R$ state independently by either $-2$ $-1$, $0$, $1$, or $2$ children with equal probability. Over 5 years, this results in 19 perturbations, not including the initial condition.  In reality, child attendance in schools will naturally fluctuate over time.  These perturbations also help in identifying the correct model by perturbing the system off of the attractor. 
%We perform inference to discover the equations for $S$ and $I$, and do not include time-series of $R$ in data matrix $\mathbf{X}_T$ or derivative matrix $\mathbf{\dot{X}}_T$ for this problem.  
%
In this example, the training and validation sets rely solely on $S$ and $I$ such that $ \mathbf{Y}_T = [\mathbf{S}(t_i)^T \; \mathbf{I}(t_i)^T]$.
The validation time-series, $\mathbf{Y}_V$ are constructed with the same number of temporal samples from a new initial condition $S_0 =15$ $I_0 =10$, $R_0= 975$. 

%Without any additional perturbations to the system, the dynamics will settle into periodic behavior, where the size of spring and fall outbreaks are the same every year. Such behavior is not only physically and biologically un-realistic, but also makes model identification impossible. The transient dynamics of the system are effectively only being sampled during the first year, after which the same dynamics are sampled over and over. To produce a more realistic time-series that better samples the possible dynamics of the system we add random `kicks' to the system at the start of each session. These perturbations take the form of adding or removing one or two children from each population $S$, $I$, and $R$. This is in addition to the constant `birth' and `death' rates defined in Eqns. \ref{eq:S}-\ref{eq:R}. For proof of concept, we do not add any noise to this system.

\subsubsection{Hybrid-SINDy discovers the switching from school breaks}

The relatively low transmission rate when school is out of session leads to an increase in the susceptible population.  As school starts, the increase in mixing between children initiates a rapid increase in the infected population, illustrated in Fig.~\ref{fig:SIR}(1\&2).  The training data for Hybrid-SINDy includes the $S$ and $I$ time-series illustrated in Fig.~\ref{fig:SIR}(6). The validation time-series is used to calculate the AIC$_c$ values.  Here, we cluster the measurement data using the coordinates $S$ and $I$, with $K = 30$ points per cluster. We use a model library containing polynomials up to 3rd order in terms of $\mathbf{X}_T$.  Two models appear with high frequency across a majority of the clusters.  The highest frequency model identifies the correct dynamical terms described in Eq. \eqref{eq:R}.  The other frequently identified model is a system with zero dynamics.

% followed by a spike in the infected population at the start of the school session as shown in the simulated time-series for the system over one year (section 2 of Fig. \ref{fig:SIR}). The full training $S$ and $I$ time-series used for model inference is shown as a function time in section 6 \ref{fig:SIR}. In addition to this training set, we simulate a validation time-series of the same length and characteristics from a new initial condition. This validation time-series is used to calculate the AIC$_c$ values, as described in the methods section.  We cluster in the data-driven coordinates $S$ vs $I$ (section 3), using $k = 30$ points per cluster, and find there are two models that appear with high frequency across clusters as the lowest AIC$_c$ model within a single cluster (section 4). The highest frequency model has the correct structure. 

Examining the coefficients for the highest frequency model over time, we identify three reoccurring sets of coefficients, illustrated in Fig.~\ref{fig:SIR}(5). The first set of recovered values correctly matches the coefficients for Eqs. \eqref{eq:S} and \eqref{eq:I} when school is out, the second set correctly recover coefficients for when school is in session, and the third are incorrect. Only the coefficient on the nonlinear transmission term, $I S$, changes value between the recovered in-school (pink) and out-of-school (maroon) transmission rates. The other coefficients (purple) are constant across the first two sets of coefficients. 

%When we plot the coefficients for the highest frequency model over time, we find three reoccurring sets of coefficients. The first set of recovered values correctly match the coefficients for Eqns \ref{eq:S} and \ref{eq:I} when school is out, the second set correctly recover coefficients for when school is in session, and the third are erroneous. Only the coefficient on the transmission term, $I S$, changes value between the recovered in-school (pink) and out-of-school (maroon) transmission rates. The other coefficients (purple) are constant across the first two sets of coefficients. 

The third set of coefficients (grey) are incorrect.  However, during these periods of time the most frequently appearing model no longer has the lowest AIC$_c$. 
The second highest frequency model $\dot{S} = 0$, $\dot{I} = 0$ has the lowest AIC$_c$ values at those times.  Additionally, the AIC$_c$ values are 4 orders of magnitude larger than those calculated for the correct model with correct coefficients. Notably, the second highest frequency model is identified by Hybrid-SINDy for regions where $S$ and $I$ are not changing because the system has reached a temporary equilibrium. This model is locally accurate, but cannot predict the validation data once a new outbreak occurs, and thus has a high (AIC$_c$ $\approx 10^{-3}$ to $1$) compared to the correct model (AIC$_c$ $\approx 10^{-6}$ to $10^{-8}$).

%The third set of coefficients (grey) are incorrect.  However, during these periods of time the most frequently appearing model no longer has the lowest AIC$_c$. 
%The second highest frequency model $dS = 0$, $dI = 0$ has the lowest AIC$_c$ values at those times, and the values are 4 orders of magnitude larger than those calculated for the correct model with correct coefficients (Data not shown). Notably, the second highest frequency model is identified by SINDy for regions where $S$ and $I$ are not changing because the system has reached a temporary equilibrium. This model is in some sense locally correct, but cannot predict the validation data once a new outbreak occurs, and thus has a high AIC$_c$ $=10^{-3}$ to $1$ compared to the correct model AIC$_c$ $=10^{-6}$ to $10^{-8}$.

\subsection{Robustness of Hybrid-SINDy to noise and cluster size}

We examine the performance of Hybrid-SINDy when varying the cluster size and noise level.   The effect of cluster-size is particularly important to understand the robustness of Hybrid-SINDy. In \S\ref{s:results}\ref{ss:massspring}, Hybrid-SINDy failed to recover the model during the transition events.  This was primarily due to the inclusion of data from both the flying and hopping dynamic regions.  In this case, the size of the regions where Hybrid-SINDy is not able to identify the correct model increases with cluster size.  Alternatively, if the cluster size is too small, the SINDy regression procedure will not be able to recover the correct model from the library.  

To investigate the impact of cluster size in SINDy's success, we perform a series of numerical experiments varying the cluster size and noise level.  We generate new sets of training and testing time-series for the mass-spring hopping model. The training set consists of time-series from 100 random initial conditions normally distributed between $x_0 \in [1,1.5]$ and $v_0 \in [0,0.5]$. The validation set consists of time-series from 10 new random initial conditions normally distributed within the same range. We divide the training set into the compression and flying subsets, avoiding the switching points.  Clusters in the flying subset are constructed by picking the time-series point with maximum position value (highest flying point), and using a nearest neighbor clustering algorithm.  By increasing $K$, the size of the clusters increase.  A similar procedure is performed during the compression phase.  Cluster sizes range from $K =10$ to $14500$.

\begin{figure}[t]
	\includegraphics[width=1\textwidth]{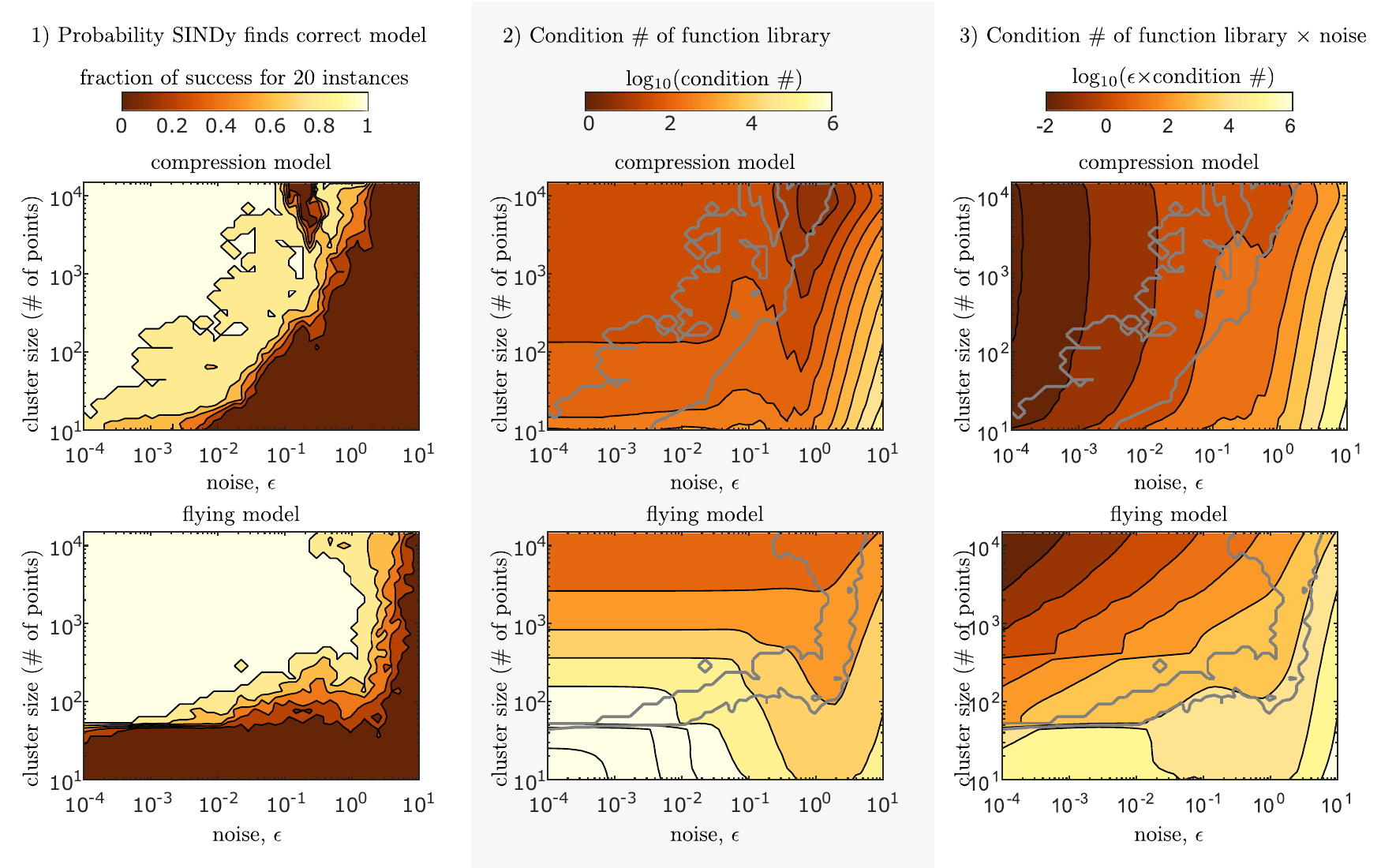}
	%%% where xxxxxx name represents "figurename.eps"
	\vspace*{-.1in}
	\caption{Success of Hybrid-SINDy on clustered compression (top) and flying (bottom) time-series points with varying noise (x-axis) and cluster size (y-axis). Sec. 1. plots shows the fraction of success in finding the correct model, over 20 noise-instances. When clusters are large and noise is low, models are recovered 100\% of the time. The color contours on Sec.2. plots indicate log$_{10}$ of the condition \# of the function library with time-series for each cluster size and noise level plugged in. Sec. 3. plots show the log$_{10}$ of the condition \# times the noise. Contours of condition \# times noise follow the contour lines for successful discovery of the model (grey). }
	\label{fig:Noise}
\end{figure}

%We generate new sets of training and testing time-series for the mass-spring hopping model. The training set consists of time-series from 100 new random initial conditions normally distributed between $x_0=$ 1 to 1.5 and $v_0=$ 0 to 0.5. The testing set consists of time-series from 10 new random initial conditions normally distributed within the same range. We divide the training set into the compression and flying subsets, so that we only select from one dynamical system at a time and avoid cross-over. Clusters in the flying subset are constructed by picking the time-series point with maximum x value (highest flying point), and using k-nearest neighbor search with increasing $k$ to find larger clusters. These clusters include increasing amounts of time-series points near "take-off" and "landing" with increasing $k$. Similarly, clusters in the compression subset we find the minimum x value (most compressed point), and use the k-nearest neighbor search to find increasingly large clusters, including points closer and closer to the transition point $x=1$. Cluster sizes range from $k =10$ to $14500$.

We also evaluated the recovery of these clusters by increasing measurement noise.  Normally-distributed noise with mean zero and $\epsilon$ from $10^{-4}$ to $10$ was added to the position, $x$, and velocity, $v$, training and testing time-series data in $\mathbf{X}$. We computed the derivatives in $\dot{\mathbf{X}}$ exactly, isolating measurements noise from the challenge of computing derivatives from noisy data.  For each cluster size and noise level, we generated 20 different noise realizations. SINDy is applied to each realization separately, and the fraction of successful model identifications are shown as the color intensity in Fig. \ref{fig:Noise}(1). With high fidelity, Hybrid-SINDy recovers the correct models for both the compression and flying clusters, when noise is relatively low and the cluster-size is relatively large, Fig. \ref{fig:Noise}(1). 
Interestingly, the cluster and noise-threshold are not the same for the compression and flying model. Recovery of the compression model varies with both noise and cluster size (the noise threshold increases for larger clusters). The cluster threshold, near $K>50$ points, and noise threshold, near $\epsilon<1$, for recovery of the flying model are independent. Notably, the flying model, which is simpler than the compression model, requires larger cluster size at the low noise limit. 

%We evaluated the recovery for these clusters with varying noise conditions. Normally-distributed noise with mean zero and $\epsilon$ from $10^{-4}$ to $10$ was added to the position, $x$ and velocity $v$ training and testing time-series sets. The derivatives were calculated exactly, not from noisy time-series. For each cluster size and noise level, we generated 20 noise instances and ran SINDy. As would be expected, SINDy recovers the correct models for both the compression and flying clusters 100\% of the time, when noise is relatively low and the cluster-size is relatively large (Sec. 1 of Fig. \ref{fig:Noise}). Perhaps less intuitively, the cluster-threshold and noise-threshold are not the same for the compression and flying model. Recovery of the compression model varies with both noise and cluster size (the noise threshold increases for larger clusters). The cluster threshold, near $k>50$ points, and noise threshold, near $\epsilon<1$, for recovery of the flying model seem to be fairly independent. Notably, the flying model, which is simpler than the compression model, requires larger cluster size at the low noise limit. 

\subsection{Condition number and noise magnitude offers insight}

To investigate the recovery patterns and the discrepancy between the compression and flying model, we calculate the condition number of $\boldsymbol{\Theta}(\mathbf{X})$ for each cluster-size and noise magnitude, as shown in Fig. \ref{fig:Noise}(2). The range of condition numbers between the two dynamical regimes are notably different.  Further, the threshold for recovery (grey) does not follow the contours for the condition number. If we instead plot contours of condition number times noise magnitude, $\kappa \epsilon$, as shown in Fig. \ref{fig:Noise}(3), the contours for successful model discovery match well. The threshold, $\kappa \epsilon$, required for discovery of the compression model is much lower than that for the flying model.

The $\kappa \epsilon$ diagnostic can be related to the noise-induced error in the least squares solution-- that is, the error in the solution of \eqref{eq:opt} with $\hat{\lambda} = 0$ and noise added to the observations $\bX$.
Because the SINDy algorithm converges to a local solution of \eqref{eq:opt} \cite{zhang2018convergence}, the closeness
of the initial least-squares iteration to the true solution gives some sense of when the algorithm will succeed. 
Let $\boldsymbol{\Xi}$ denote the true solution and $\delta \bXi^\LS$ denote the difference between the true solution and the least-squares  solution for noisy data. Then 
\begin{equation} 
\frac{\|\delta \bXi^\LS \|_2}{\|\bXi \|_2}
\leq \frac{C \kappa \epsilon}{1- C\kappa \epsilon} \; ,
\label{eq:bound}
\end{equation}
for some constant $C$ which depends only on the library functions. 
Note that the condition number, $\kappa$, depends on the sampling (cluster size) and choice of library functions.  
The complex interplay between the magnitude of noise, sampling schemes, and choice of SINDy library in (\ref{eq:bound}) provides a threshold for when we expect Hybrid-SINDy to recover the true solution.
See \S \ref{ss:appendix} for a more detailed discussion.
The plots in Fig.~\ref{fig:Noise}(3) show that this diagnostic threshold correlates well with the empirical performance of SINDy.
It remains unclear why the particular value of $\kappa \epsilon$ for which the algorithm succeeds is three orders of magnitude higher for the flying regime than that for the compression regimes. 
Intuitively, there are more terms to recover and these terms have a high contrast.
However, these considerations do not fully account for the difference in the observed behaviors.
In \S\ref{ss:appendix}, we provide some further intuition as to why the regression problem in the compression regime is more difficult than the problem in the flying regime.
A more fine-grained analysis is required, taking into account the heterogeneous effect that noise in the observations has on the values of the library functions.
%Thus, for a fixed ratio of largest and smallest nonzero
%coefficient and a fixed level of sparsity,
%the condition number times noise, $\kappa \epsilon$,
%is a natural diagnostic for the performance of SINDy.
% 

%Elevator pitch about why.

\section{Discussion and conclusions}

Characterizing the complex and dynamic interactions of physical and biological systems is essential for designing intervention and control strategies.
For example, understanding infectious disease transmission across human populations has led to better informed large-scale vaccination campaigns~\cite{upfill2014predictive,mercer2017spatial}, vector control programs~\cite{nikolov2016malaria,eckhoff2017impact}, and surveillance activities~\cite{nikolov2016malaria,gerardin2017effectiveness}.  
The increasing availability of measurement data, computational resources, and data storage capacity enables new data-driven methodologies for characterization of these systems.  
Recent methodological innovations for identifying nonlinear dynamical systems from data have been broadly successful in a wide-variety of applications including fluid dynamics~\cite{Mezic2005nd}, epidemiology~\cite{proctor2015discovering}, metabolic networks~\cite{mangan2016inferring}, and ecological systems~\cite{Sugihara1990nature,Sugihara2012}.
%
%Recent and currently popular methodologies include neural nets~\cite{goodfellow2016deep}, Koopman operator theory~\cite{Koopman1931pnas,Mezic2005nd,Williams2015jnls}, Dynamic Mode Decomposition~\cite{Rowley2009jfm,Tu2014jcd,Proctor2016siads,Kutz2016book}, and cross-convergent mapping~\cite{Sugihara2012}.  
%
The recently developed SINDy methodology identifies nonlinear models from data, offers a parsimonious and interpretable model representation~\cite{Brunton2016pnas}, and generalizes well to realistic constraints such as limited and noisy data~\cite{mangan2016inferring,Rudy2017a,Mangan2017}.
Broadly, SINDy is a data-analysis and modeling tool that can provide insight into mechanism as well as prediction.  
Despite this substantial and encouraging progress, the characterization of nonlinear systems from data is incomplete.  
Complex systems that exhibit switching between dynamical regimes have been far less studied with these methods, despite the ubiquity of these phenomena in physical, engineered, and biological systems~\cite{van2000introduction,holmes2006dynamics}.  

The primary contribution of this work is the generalization of SINDy to identify hybrid systems and their switching behavior.
We call this new methodology Hybrid-SINDy.
%
%In this era of big data, most scientific and engineering communities are rapidly innovating in this space.
%
By characterizing the similarity among data points, we identify clusters in measurement space using an unsupervised learning technique.
A set of SINDy models is produced across clusters, and the highest frequency and most informative, predictive models are selected.
We demonstrate the success of this algorithm on two modern examples of hybrid systems~\cite{holmes2006dynamics,grenfell1994measles}: the state-dependent switching of a hopping robot and the time-dependent switching of disease transmission dynamics for children in-school and on-vacation.

For the hopping robot, Hybrid-SINDy correctly identifies the flight and compression regimes.   
SINDy is able to construct candidate nonlinear models from data drawn across the entire time-series, but restricted to measurements similar in measurement space.  
This innovation allows data to be clustered based on the underlying dynamics and nonlinear geometry of trajectories, enabling the use of regression-based methods such as SINDy.
The method is also quite intuitive for state-dependent hybrid systems; phase-space is effectively partitioned based on the similarity in measurement data.
Moreover, this equation-free method is consistent with the underlying theory of hybrid dynamical systems by establishing charts where distinct nonlinear dynamical regimes exist between transition events.
We also demonstrate that Hybrid-SINDy correctly identifies time-dependent hybrid systems from a subset of all of the phase variables.
We can identify the SIR system with separate transmission rates among children during in-school versus on-vacation mixing patterns, based solely on the susceptible and infected measurements of the system.  
For both examples, we show that the model error characteristics and the library of candidate models help illustrate the switching behavior even in the presence of additive measurement noise.   
These examples illustrate the adaptability of the method to realistic measurements and complex system behaviors.

Hybrid-SINDy incorporates the fundamental elements of a broad number of other methodologies.
The method builds a library of features from measurement data to better predict the future measurement.  
Variations of this augmentation process has been widely explored over the last few decades, notably in the control theoretic community with delay embeddings~\cite{kalman:1965,Juang1985jgcd,Phan1992jas,Katayama.2005}, Carleman linearization~\cite{kowalski1991nonlinear}, and nonlinear autoregressive models~\cite{box2015time}.  
More recently, machine-learning and computer science approaches often refer to the procedure as feature engineering.      
Constraining the input data is another well-known approach to identify more informative and predictive models.  
Examples include windowing the data in time for autoregressive moving average models or identifying similarity among measurements based on the Takens' embedding theorem for delay embeddings of chaotic dynamical systems~\cite{Takens1981lnm,Sugihara1994ptrsla,Sugihara2012} .
With Hybrid-SINDy, we integrate and adapt a number of these components to construct an algorithm that can indentify nonlinear dynamical systems and switching between dynamical regimes. 

There are limitations and challenges to the wide-spread adoption of our method.
The method is fundamentally data-driven, requiring an adequate amount of data for each dynamical regime to perform the SINDy regression.
We also rely on having access to a sufficient number of measurement variables to construct the nonlinear dynamics, even with the inclusion of delay embeddings. 
These measurements also need to be in a coordinate frame to allow for a parsimonious description of the dynamics.
In order to test the robustness of our results, we evaluate the condition number and noise magnitude as a numerical diagnostic for evaluating the output of Hybrid-SINDy.
However, despite developing a rigorous mathematical connection between this diagnostic and numerically solving the SINDy regression, we discovered that there does not exist a specific threshold number that generalizes across models and library choice.

Despite these limitations, Hybrid-SINDy is a novel step toward a general method for identifying hybrid nonlinear dynamical systems from data.  
We have mitigated a number of the numerical challenges by incoporating information theoretic criteria to manage uncertainty and offering a procedure to validate the results against cluster size and noise magnitude.
Looking ahead, discovering a general criteria that holds across a wide-variety of applictions and models will be essential for the wide-spread adoption of this methodology. 
Further, we foresee the innovative work around data-driven identification of nonlinear manifolds as another important research direction for Hybrid-SINDy~\cite{COIFMAN20065,Giannakis2012pnas}.
%
% We believe Hybrid-SINDy offers a new perspective for evaluating data from systems that switch dynamic behavior and provides intuition about these complex systems.

\section{Appendix A: Bound derivation}
\label{ss:appendix}
Zhang and Schaeffer \cite[Theorem 2.5]{zhang2018convergence} showed that the SINDy hard-thresholding procedure converges to a local solution of \eqref{eq:opt} with $R(\cdot) = \|\cdot\|_0$.
Because that problem is nonconvex, a local solution may or may not be equal to the true global solution.
We are interested in characterizing when the initial guess for SINDy is ``close'' to the exact sparse solution.

For each value of $\hat{\lambda}$, we initialize SINDy with the least squares solution, i.e. the solution of \eqref{eq:opt} with $\hat{\lambda}=0$.
Noise is added to the observations $\bX$ alone and $\dot{\bX}$ is without noise.
Let $\bX + \delta\bX$ denote the noisy data and let $\delta \bTheta \coloneqq \bTheta( \bX + \delta \bX) - \bTheta(\bX)$ denote the perturbation in the resulting library. 
For the sake of simplicity, we will assume that $ \| \delta \bTheta \|_2/\| \bTheta \|_2 \leq C \| \delta \bX \|_2 / \| \bX \|_2  \leq C \epsilon$, where $\epsilon$ is the noise level and $C$ depends only on the choice of library functions. 
We further assume that $\bTheta$ and
$\bTheta + \delta \bTheta$ are full rank and that
$\dot{\bX} = \bTheta \bXi$, i.e. that $\dot{\bX}$ is
in the range of $\bTheta$ so that the true solution
$\bXi$ satisfies $\bXi = \bTheta^\dagger \dot{\bX}$, where
$\dagger$ denotes the Moore-Penrose pseudo-inverse.

The solution of the noisy least-squares problem is then
$\bXi + \delta \bXi^\LS = \left ( \bTheta + \delta \bTheta \right )^\dagger
\dot{\bX}$ ,
where $\delta \bXi^\LS$ denotes the resulting error.
Let $\kappa = \| \bTheta \|_2 \| \bTheta^\dagger \|_2$ denote the
condition number of $\bTheta$. We have the bound \eqref{eq:bound} from the main text,

\begin{equation} \label{eq:l2bound}
\frac{\|\delta \bXi^\LS \|_2}{\|\bXi \|_2}
\leq \frac{C\kappa \epsilon}{1-C\kappa \epsilon} \; ,
\end{equation}
provided that $C\kappa \epsilon < 1$. The derivation of \eqref{eq:l2bound}
is nontrivial; for a reference, see
\cite[Theorem 5.1]{wedin1973perturbation}.

To see why the flying model is easier to recover than the compression model for a given value of $\kappa \epsilon$,we consider a single step of hard thresholding.  
For this derivation, we consider $\bXi$ to be a vector; this assumption holds when $\bXi$ is not a vector, since we typically consider solving for each column of $\bXi$ independently in the SINDy regression.
Let $c$ denote the size of the smallest nonzero coefficient in $\bXi$,
i.e. $c = \min_{i,j \textrm{ s.t. } \Xi_{ij} \neq 0} \left |\Xi_{ij} \right |$.
A single step of hard thresholding will succeed in finding the true
support of $\bXi$ using the threshold $c/2$ when $\| \delta \bXi^\LS \|_\infty$
is smaller than $c/2$. Let $k$ be the number of non-zero entries in
$\bXi$. Observing that $\|\bXi\|_2 \leq \sqrt{k}\|\bXi\|_\infty$,
we have

\begin{equation}
\frac{\| \delta \bXi^\LS \|_\infty}{c}
\leq \frac{\sqrt{k} \|\bXi\|_\infty}{c}
\frac{\|\delta \bXi^\LS \|_2}{\|\bXi \|_2}
\leq \frac{\sqrt{k} \|\bXi\|_\infty}{c}
\frac{C\kappa \epsilon}{1-C\kappa \epsilon} \; .
\end{equation}

We see that the number of nonzero coefficients and the ratio of the largest to smallest coefficients in the true solution affect the success of a single step of hard thresholding.
Intuitively, then, the compression model is more difficult to recover than the simpler flying model.
However, the factor $\sqrt{k} \|\bXi\|_\infty/c$ only accounts for about an order of magnitude of the discrepancy in the $\kappa\epsilon$ threshold at which SINDy correctly recovered compression and flying models in Fig. \ref{fig:Noise}.
A likely culprit for the remaining difference is the variation in the effect of noise on different basis functions in the library.
For example, adding noise to $\bX$ has no effect on the constant term,  but will be magnified by a quadratic term.
A more fine-grained analysis of the error corresponding to the specific functions in the model could account for the remaining discrepancy.

%\ack{Insert acknowledgment text here.}

\section*{Acknowledgements}
{JLP and NMM would like to thank Bill and Melinda Gates for their active support of the Institute for Disease Modeling and their sponsorship through the Global Good Fund.  JNK and TA acknowledge support from the Air Force Office of Scientific Research (FA9550-15-1-0385,  FA9550-17-1-0329). SLB and JNK acknowledge support from the Defense Advanced Research Projects Agency (DARPA contract HR0011-16-C-0016). SLB acknowledges support from the Army Research Office (W911NF-17-1-0422) }

%%%%%%%%%% Insert bibliography here %%%%%%%%%%%%%%
\bibliographystyle{prsb}
\bibliography{references}

\begin{thebibliography}{10}
\expandafter\ifx\csname urlstyle\endcsname\relax
  \providecommand{\doi}[1]{doi:\discretionary{}{}{}#1}\else
  \providecommand{\doi}{doi:\discretionary{}{}{}\begingroup
  \urlstyle{rm}\Url}\fi

\bibitem{keeling2001seasonally}
Keeling, M.~J., Rohani, P. \& Grenfell, B.~T., 2001 Seasonally forced disease
  dynamics explored as switching between attractors.
\newblock \emph{Physica D: Nonlinear Phenomena} \textbf{148}, 317--335.

\bibitem{holmes2006dynamics}
Holmes, P., Full, R.~J., Koditschek, D. \& Guckenheimer, J., 2006 The dynamics
  of legged locomotion: Models, analyses, and challenges.
\newblock \emph{SIAM review} \textbf{48}, 207--304.

\bibitem{dobson2007complex}
Dobson, I., Carreras, B.~A., Lynch, V.~E. \& Newman, D.~E., 2007 Complex
  systems analysis of series of blackouts: Cascading failure, critical points,
  and self-organization.
\newblock \emph{Chaos: An Interdisciplinary Journal of Nonlinear Science}
  \textbf{17}, 026103.

\bibitem{li2014communication}
Li, H., Dimitrovski, A.~D., Song, J.~B., Han, Z. \& Qian, L., 2014
  Communication infrastructure design in cyber physical systems with
  applications in smart grids: A hybrid system framework.
\newblock \emph{IEEE Communications Surveys and Tutorials} \textbf{16},
  1689--1708.

\bibitem{Akaike1973}
Akaike, H., 1973 {Information theory and an extension of the maximum likelihood
  principle}.
\newblock In \emph{Petrov, B.N.; Cs{\'{a}}ki, F., 2nd International Symposium
  on Information Theory, Tsahkadsor, Armenia, USSR, September 2-8, 1971}, pp.
  267--281. Budapest: Akad{\'{e}}miai Kiad{\'{o}}.

\bibitem{Akaike1974}
Akaike, H., 1974 {A New Look at the Statistical Model Identification}.
\newblock \emph{IEEE Transactions on Automatic Control} \textbf{19}, 716--723.
\newblock ISSN 15582523.
\newblock (\doi{10.1109/TAC.1974.1100705}).

\bibitem{nakamura2006comparative}
Nakamura, T., Judd, K., Mees, A.~I. \& Small, M., 2006 A comparative study of
  information criteria for model selection.
\newblock \emph{International Journal of Bifurcation and Chaos} \textbf{16},
  2153--2175.

\bibitem{lillacci2010parameter}
Lillacci, G. \& Khammash, M., 2010 Parameter estimation and model selection in
  computational biology.
\newblock \emph{PLoS Computational Biology} \textbf{6}, e1000696.

\bibitem{penny2012comparing}
Penny, W.~D., 2012 Comparing dynamic causal models using {AIC}, {BIC} and free
  energy.
\newblock \emph{Neuroimage} \textbf{59}, 319--330.

\bibitem{Brunton2016pnas}
Brunton, S.~L., Proctor, J.~L. \& Kutz, J.~N., 2016 Discovering governing
  equations from data by sparse identification of nonlinear dynamical systems.
\newblock \emph{Proceedings of the National Academy of Sciences} \textbf{113},
  3932--3937.

\bibitem{Mangan2017}
Mangan, N.~M., Kutz, J.~N., Brunton, S.~L. \& Proctor, J.~L., 2017 {Model
  selection for dynamical systems via sparse regression and information
  criteria}.
\newblock \emph{Proceedings of the Royal Society of London A: Mathematical,
  Physical and Engineering Sciences} \textbf{473}.

\bibitem{mangan2016inferring}
Mangan, N.~M., Brunton, S.~L., Proctor, J.~L. \& Kutz, J.~N., 2016 Inferring
  biological networks by sparse identification of nonlinear dynamics.
\newblock \emph{IEEE Transactions on Molecular, Biological and Multi-Scale
  Communications} \textbf{2}, 52--63.

\bibitem{Rudy2017a}
Rudy, S.~H., Brunton, S.~L., Proctor, J.~L. \& Kutz, J.~N., 2017 {Data-driven
  discovery of partial differential equations}.
\newblock \emph{Science Advances} \textbf{3}.
\newblock ISSN 23752548.
\newblock (\doi{10.1126/sciadv.1602614}).

\bibitem{kalman:1965}
Ho, B.~L. \& Kalman, R.~E., 1965 Effective construction of linear
  state-variable models from input/output data.
\newblock In \emph{Proceedings of the 3rd Annual Allerton Conference on Circuit
  and System Theory}, pp. 449--459.

\bibitem{Juang1985jgcd}
Juang, J.~N. \& Pappa, R.~S., 1985 {An eigensystem realization algorithm for
  modal parameter identification and model reduction}.
\newblock \emph{Journal of Guidance, Control, and Dynamics} \textbf{8},
  620--627.

\bibitem{Phan1992jas}
Phan, M., Juang, J.~N. \& Longman, R.~W., 1992 Identification of
  linear-multivariable systems by identification of observers with assigned
  real eigenvalues.
\newblock \emph{The Journal of the Astronautical Sciences} \textbf{40},
  261--279.

\bibitem{Katayama.2005}
T., K., 2005 \emph{{Subspace Methods for System Identification}}.
\newblock Springer-Verlag London.

\bibitem{Sugihara1990nature}
Sugihara, G. \& May, R.~M., 1990 Nonlinear forecasting as a way of
  distinguishing chaos from measurement error in time series.
\newblock \emph{Nature} \textbf{344}, 734--741.

\bibitem{Sugihara1994ptrsla}
Sugihara, G., 1994 Nonlinear forecasting for the classification of natural time
  series.
\newblock \emph{Philosophical Transactions of the Royal Society of London A}
  \textbf{348}, 477--495.

\bibitem{brunton2017chaos}
Brunton, S.~L., Brunton, B.~W., Proctor, J.~L., Kaiser, E. \& Kutz, J.~N., 2017
  Chaos as an intermittently forced linear system.
\newblock \emph{Nature Communications} \textbf{8}, 19.

\bibitem{Sugihara2012}
Sugihara, G., May, R., Ye, H., Hsieh, C., Deyle, E., Fogarty, M. \& Munch, S.,
  2012 Detecting causality in complex ecosystems.
\newblock \emph{Science} \textbf{338}, 496--500.

\bibitem{kowalski1991nonlinear}
Kowalski, K. \& Steeb, W.-H., 1991 \emph{Nonlinear dynamical systems and
  Carleman linearization}.
\newblock World Scientific.

\bibitem{Schmid2010jfm}
Schmid, P.~J., 2010 Dynamic mode decomposition of numerical and experimental
  data.
\newblock \emph{Journal of Fluid Mechanics} \textbf{656}, 5--28.
\newblock ISSN 0022-1120.

\bibitem{Rowley2009jfm}
Rowley, C., Mezi{\'{c}}, I., Bagheri, S., Schlatter, P. \& Henningson, D.~S.,
  2009 {Spectral analysis of nonlinear flows}.
\newblock \emph{J. Fluid Mech.} \textbf{645}, 115--127.

\bibitem{kutz2016dynamic}
Kutz, J.~N., Brunton, S.~L., Brunton, B.~W. \& Proctor, J.~L., 2016
  \emph{Dynamic mode decomposition: data-driven modeling of complex systems},
  volume 149.
\newblock SIAM.

\bibitem{Mezic2005nd}
Mezi{\'c}, I., 2005 Spectral properties of dynamical systems, model reduction
  and decompositions.
\newblock \emph{Nonlinear Dynamics} \textbf{41}, 309--325.

\bibitem{Mezic2013arfm}
Mezic, I., 2013 Analysis of fluid flows via spectral properties of the
  {K}oopman operator.
\newblock \emph{Annual Review of Fluid Mechanics} \textbf{45}, 357--378.

\bibitem{COIFMAN20065}
Coifman, R.~R. \& Lafon, S., 2006 Diffusion maps.
\newblock \emph{Applied and Computational Harmonic Analysis} \textbf{21}, 5 --
  30.
\newblock ISSN 1063-5203.
\newblock (\doi{10.1016/j.acha.2006.04.006}).
\newblock Special Issue: Diffusion Maps and Wavelets.

\bibitem{Giannakis2012pnas}
Giannakis, D. \& Majda, A.~J., 2012 Nonlinear {L}aplacian spectral analysis for
  time series with intermittency and low-frequency variability.
\newblock \emph{Proceedings of the National Academy of Sciences} \textbf{109},
  2222--2227.

\bibitem{box2015time}
Box, G.~E., Jenkins, G.~M., Reinsel, G.~C. \& Ljung, G.~M., 2015 \emph{Time
  series analysis: forecasting and control}.
\newblock John Wiley \& Sons.

\bibitem{macevsic2017koopman}
Ma{\'c}e{\v{s}}i{\'c}, S., {\v{C}}rnjari{\'c}-{\v{Z}}ic, N. \& Mezi{\'c}, I.,
  2017 Koopman operator family spectrum for nonautonomous systems-part 1.
\newblock \emph{arXiv preprint arXiv:1703.07324} .

\bibitem{Linderman2016}
Linderman, S.~W., Miller, A.~C., Adams, R.~P., Blei, D.~M., Paninski, L. \&
  Johnson, M.~J., 2016 {Recurrent switching linear dynamical systems}.
\newblock \emph{arXiv preprint arXiv:1610.08466} .

\bibitem{yair2017reconstruction}
Yair, O., Talmon, R., Coifman, R.~R. \& Kevrekidis, I.~G., 2017 Reconstruction
  of normal forms by learning informed observation geometries from data.
\newblock \emph{Proceedings of the National Academy of Sciences} \textbf{114},
  E7865--E7874.

\bibitem{back1993dynamical}
Back, A., Guckenheimer, J. \& Myers, M., 1993 A dynamical simulation facility
  for hybrid systems.
\newblock In \emph{Hybrid systems}, pp. 255--267. Springer.

\bibitem{van2000introduction}
Van Der~Schaft, A.~J. \& Schumacher, J.~M., 2000 \emph{An introduction to
  hybrid dynamical systems}, volume 251.
\newblock Springer London.

\bibitem{Hesterberg2008}
Hesterberg, T., Choi, N.~H., Meier, L. \& Fraley, C., 2008 {Least angle and
  $\ell_1$ penalized regression: A review}.
\newblock \emph{Statistics Surveys} \textbf{2}, 61--93.
\newblock ISSN 1935-7516.
\newblock (\doi{10.1214/08-SS035}).

\bibitem{Schaeffer2017a}
Schaeffer, H., 2017 {Learning partial differential equations via data discovery
  and sparse optimization}.
\newblock \emph{Proceedings of the Royal Society A: Mathematical, Physical and
  Engineering Science} \textbf{473}, 20160446.
\newblock ISSN 1364-5021.
\newblock (\doi{10.1098/rspa.2016.0446}).

\bibitem{Tran2016arxiv}
TRAN, G. \& WARD, R., 2017 {EXACT RECOVERY OF CHAOTIC SYSTEMS FROM HIGHLY
  CORRUPTED DATA}.
\newblock \emph{Multiscale Modeling {\&} Simulation} \textbf{15}, 1108--1129.

\bibitem{Schaeffer2017b}
Schaeffer, H. \& McCalla, S.~G., 2017 {Sparse model selection via integral
  terms}.
\newblock \emph{Physical Review E} \textbf{96}, 023302.
\newblock ISSN 24700053.
\newblock (\doi{10.1103/PhysRevE.96.023302}).

\bibitem{Pantazis2017}
Pantazis, Y. \& Tsamardinos, I., 2017 A unified approach for sparse dynamical
  system inference from temporal measurements.
\newblock \emph{arXiv preprint arXiv:1710.00718} .

\bibitem{Schaeffer2017}
Schaeffer, H., Tran, G. \& Ward, R., 2017 Learning dynamical systems and
  bifurcation via group sparsity.
\newblock \emph{arXiv preprint arXiv:1611.03271} pp. 1--16.

\bibitem{Rudy2018arxiv}
Rudy, S., Alla, A., Brunton, S.~L. \& Kutz, J.~N., 2018 Data-driven
  identification of parametric partial differential equations.
\newblock \emph{arXiv preprint arXiv:1806.00732} .

\bibitem{Loiseau2017jfm}
Loiseau, J.-C. \& Brunton, S.~L., 2018 Constrained sparse {Galerkin}
  regression.
\newblock \emph{Journal of Fluid Mechanics} \textbf{838}, 42--67.

\bibitem{burnham2002}
Burnham, K. \& Anderson, D., 2002 \emph{{Model Selection and Multi-Model
  Inference}}.
\newblock Springer, 2nd edition.

\bibitem{Kuepfer2007}
Kuepfer, L., Peter, M., Sauer, U. \& Stelling, J., 2007 {Ensemble modeling for
  analysis of cell signaling dynamics.}
\newblock \emph{Nature Biotechnology} \textbf{25}, 1001--1006.
\newblock ISSN 1087-0156.
\newblock (\doi{10.1038/nbt1330}).

\bibitem{Hjorth2008}
Claeskens, G. \& Hjorth, N.~L., 2008 \emph{{Model Selection and Model
  Averaging}}.
\newblock Cambridge University Press.

\bibitem{Schaber2011}
Schaber, J., Fl{\"{o}}ttmann, M., Li, J., Tiger, C.~F., Hohmann, S. \& Klipp,
  E., 2011 {Automated ensemble modeling with modelMaGe: Analyzing feedback
  mechanisms in the Sho1 branch of the HOG pathway}.
\newblock \emph{PLoS ONE} \textbf{6}, 1--7.
\newblock ISSN 19326203.
\newblock (\doi{10.1371/journal.pone.0014791}).

\bibitem{woodward2004epidemiology}
Woodward, M., 2004 \emph{{Epidemiology: Study Design and Data Analysis, Second
  Edition}}.
\newblock Chapman {\&} Hall/CRC Texts in Statistical Science. Taylor {\&}
  Francis.
\newblock ISBN 9781584884156.

\bibitem{Blake22072014}
Blake, I.~M., Martin, R., Goel, A., Khetsuriani, N., Everts, J., Wolff, C.,
  Wassilak, S., Aylward, R.~B. \& Grassly, N.~C., 2014 {The role of older
  children and adults in wild poliovirus transmission}.
\newblock \emph{Proceedings of the National Academy of Sciences} \textbf{111},
  10604--10609.
\newblock ISSN 0027-8424, 1091-6490.
\newblock (\doi{10.1073/pnas.1323688111}).

\bibitem{Schmidt2009science}
Schmidt, M. \& Lipson, H., 2009 Distilling free-form natural laws from
  experimental data.
\newblock \emph{Science} \textbf{324}, 81--85.

\bibitem{Buchel2013}
B{\"{u}}chel, F., Rodriguez, N., Swainston, N., Wrzodek, C., Czauderna, T.,
  Keller, R., Mittag, F., Schubert, M., Glont, M., Golebiewski, M.
  \emph{et~al.}, 2013 {Path2Models: large-scale generation of computational
  models from biochemical pathway maps.}
\newblock \emph{BMC Systems Biology} \textbf{7}, 116.
\newblock ISSN 1752-0509.
\newblock (\doi{10.1186/1752-0509-7-116}).

\bibitem{Cohen2015}
Cohen, P.~R., 2015 {DARPA's Big Mechanism program.}
\newblock \emph{Physical Biology} \textbf{12}, 045008.
\newblock ISSN 1478-3975.
\newblock (\doi{10.1088/1478-3975/12/4/045008}).

\bibitem{Kullback1951}
Kullback, S. \& Leibler, R.~A., 1951 {On Information and Sufficiency}.
\newblock \emph{The Annals of Mathematical Statistics} \textbf{22}, 79--86.
\newblock ISSN 0003-4851.
\newblock (\doi{10.1214/aoms/1177729694}).

\bibitem{schwarz1978estimating}
Schwarz, G., 1978 Estimating the dimension of a model.
\newblock \emph{The Annals of Statistics} \textbf{6}, 461--464.

\bibitem{bishop2006pattern}
Bishop, C.~M. \& Others, 2006 \emph{{Pattern recognition and machine
  learning}}, volume~1.
\newblock Springer New York.

\bibitem{linde2005dic}
Linde, A., 2005 {DIC} in variable selection.
\newblock \emph{Statistica Neerlandica} \textbf{59}, 45--56.

\bibitem{rissanen1978modeling}
Rissanen, J., 1978 Modeling by shortest data description.
\newblock \emph{Automatica} \textbf{14}, 465--471.

\bibitem{killick2012}
Killick, R., Fearnhead, P. \& Eckley, I.~A., 2012 {Optimal Detection of
  Changepoints With a Linear Computational Cost}.
\newblock \emph{Journal of the American Statistical Association} \textbf{107},
  1590--1598.
\newblock ISSN 0162-1459.
\newblock (\doi{10.1080/01621459.2012.737745}).

\bibitem{ruina1998nonholonomic}
Ruina, A., 1998 Nonholonomic stability aspects of piecewise holonomic systems.
\newblock \emph{Reports on Mathematical Physics} \textbf{42}, 91--100.

\bibitem{full1999templates}
Full, R.~J. \& Koditschek, D.~E., 1999 Templates and anchors: neuromechanical
  hypotheses of legged locomotion on land.
\newblock \emph{Journal of Experimental Biology} \textbf{202}, 3325--3332.

\bibitem{saranli2001rhex}
Saranli, U., Buehler, M. \& Koditschek, D.~E., 2001 Rhex: A simple and highly
  mobile hexapod robot.
\newblock \emph{The International Journal of Robotics Research} \textbf{20},
  616--631.

\bibitem{raibert2008bigdog}
Raibert, M., Blankespoor, K., Nelson, G. \& Playter, R., 2008 Bigdog, the
  rough-terrain quadruped robot.
\newblock \emph{IFAC Proceedings Volumes} \textbf{41}, 10822--10825.

\bibitem{grenfell1994measles}
Grenfell, B., Kleckzkowski, A., Ellner, S. \& Bolker, B., 1994 Measles as a
  case study in nonlinear forecasting and chaos.
\newblock \emph{Phil. Trans. R. Soc. Lond. A} \textbf{348}, 515--530.

\bibitem{broadfootmeasles}
Broadfoot, K. \& Keeling, M. Measles epidemics in vaccinated populations .

\bibitem{zhang2018convergence}
Zhang, L. \& Schaeffer, H., 2018 On the convergence of the {SIND}y algorithm.
\newblock \emph{arXiv preprint arXiv:1805.06445} .

\bibitem{upfill2014predictive}
Upfill-Brown, A.~M., Lyons, H.~M., Pate, M.~A., Shuaib, F., Baig, S., Hu, H.,
  Eckhoff, P.~A. \& Chabot-Couture, G., 2014 Predictive spatial risk model of
  poliovirus to aid prioritization and hasten eradication in {N}igeria.
\newblock \emph{BMC Medicine} \textbf{12}, 92.

\bibitem{mercer2017spatial}
Mercer, L.~D., Safdar, R.~M., Ahmed, J., Mahamud, A., Khan, M.~M., Gerber, S.,
  O'Leary, A., Ryan, M., Salet, F., Kroiss, S.~J. \emph{et~al.}, 2017 Spatial
  model for risk prediction and sub-national prioritization to aid poliovirus
  eradication in {P}akistan.
\newblock \emph{BMC Medicine} \textbf{15}, 180.

\bibitem{nikolov2016malaria}
Nikolov, M., Bever, C.~A., Upfill-Brown, A., Hamainza, B., Miller, J.~M.,
  Eckhoff, P.~A., Wenger, E.~A. \& Gerardin, J., 2016 Malaria elimination
  campaigns in the {L}ake {K}ariba region of {Z}ambia: a spatial dynamical
  model.
\newblock \emph{PLoS Computational Biology} \textbf{12}, e1005192.

\bibitem{eckhoff2017impact}
Eckhoff, P.~A., Wenger, E.~A., Godfray, H. C.~J. \& Burt, A., 2017 Impact of
  mosquito gene drive on malaria elimination in a computational model with
  explicit spatial and temporal dynamics.
\newblock \emph{Proceedings of the National Academy of Sciences} \textbf{114},
  E255--E264.

\bibitem{gerardin2017effectiveness}
Gerardin, J., Bever, C.~A., Bridenbecker, D., Hamainza, B., Silumbe, K.,
  Miller, J.~M., Eisele, T.~P., Eckhoff, P.~A. \& Wenger, E.~A., 2017
  Effectiveness of reactive case detection for malaria elimination in three
  archetypical transmission settings: a modelling study.
\newblock \emph{Malaria journal} \textbf{16}, 248.

\bibitem{proctor2015discovering}
Proctor, J.~L. \& Eckhoff, P.~A., 2015 Discovering dynamic patterns from
  infectious disease data using dynamic mode decomposition.
\newblock \emph{International health} \textbf{7}, 139--145.

\bibitem{Takens1981lnm}
Takens, F., 1981 Detecting strange attractors in turbulence.
\newblock \emph{Lecture Notes in Mathematics} \textbf{898}, 366--381.

\bibitem{wedin1973perturbation}
Wedin, P.-{\AA}., 1973 Perturbation theory for pseudo-inverses.
\newblock \emph{BIT Numerical Mathematics} \textbf{13}, 217--232.

\end{thebibliography}

\end{document}